\documentclass[11pt]{article}
\usepackage{latexsym}
\usepackage{amsmath}

\usepackage{graphicx}
\usepackage{amsmath}
\usepackage{amssymb}

\def\N{\mathbb{N}}

\def\P{\mathbb{P}}
\newcommand{\sign}{\mathrm{sgn}}
\def\n{^{(n)}}
\def\np1{^{(n+1)}}

\def\njp1{^{(n_j+1)}}
\def\j{^{(j)}}
\def\S{\mathbb{S}}
\def\nell{n'_{\ell}}
\def\hell{h'^{(\ell)}}
\def\uell{u'^{(\ell)}}
\def\tuell{\tilde{u}'^{(\ell)}}
\def\vell{v'^{(\ell)}}

\newcommand{\wlim}{\mbox{w}\!\mbox{-}\!\mbox{lim}}

\newtheorem{Theorem}{Theorem}[section]
\newtheorem{Definition}[Theorem]{Definition}
\newtheorem{Lemma}[Theorem]{Lemma}
\newtheorem{Proposition}[Theorem]{Proposition}

\newenvironment{Remark}{\stepcounter{Theorem}\noindent {\bf Remark \theTheorem\ }}{\hfill$\Box$} 

\begin{document}

\author{Ingrid Daubechies, Massimo Fornasier, and Ignace Loris}
\title{Accelerated Projected Gradient Method for Linear Inverse
Problems with Sparsity Constraints}

\maketitle

\begin{abstract}
Regularization of ill-posed linear inverse problems via $\ell_1$
penalization has been proposed for cases where the solution is known
to be (almost) sparse. One way to obtain the minimizer of such an
$\ell_1$ penalized functional is via an iterative soft-thresholding
algorithm. We propose an alternative implementation to
$\ell_1$-constraints, using a 
gradient method, with
projection on $\ell_1$-balls. The corresponding algorithm uses again
iterative soft-thresholding, now with a variable thresholding
parameter. We also propose accelerated versions of this iterative
method, using ingredients of the (linear) steepest descent method.
We prove convergence in norm for one of these projected gradient methods,
without and with acceleration.
\end{abstract}

\section{Introduction}

Our main concern in this paper is the construction of iterative
algorithms to solve inverse problems with an $\ell_1$-penalization
or an $\ell_1$-constraint, and that converge faster than the
iterative algorithm proposed in \cite{dadede04} (see also formulas
\eqref{DDD04iteration} and \eqref{soft-thresholding} below). Before we
get into technical details, we introduce here the background,
framework, and notations for our work.

In many practical problems, one cannot observe directly the
quantities of most interest; instead their values have to be
inferred from their effect on observable quantities. When this
relationship between observable $y$ and interesting quantity $f$ is
(approximately) linear, as it is in surprisingly many cases, the
situation can be modeled mathematically by the equation
\begin{equation}
y \,=\, A f~, \label{eq:1}
\end{equation}
where $A$ is a linear operator mapping a vector space $\mathcal K$
(which we assume to contain all possible ``objects'' $f$) to a
vector space $\mathcal H$ (which contains all possible data $y$).
The vector spaces $\mathcal K$ and $\mathcal H$ can be finite-- or
infinite--dimensional; in the latter case, we assume that $\mathcal
K$ and $\mathcal H$ are (separable) Hilbert spaces, and that $A:
\mathcal K \to \mathcal H$ is a bounded linear operator. Our main
goal consists in reconstructing the (unknown) element $f \in
\mathcal K$, when we are given $y$. If $A$ is a ``nice'', easily
invertible operator, and if the data $y$ are free of noise, then
this is a trivial task. Often, however, the mapping $A$ is ill-conditioned
or not
invertible. Moreover, typically (\ref{eq:1}) is
only an idealized version in which noise has been neglected; a more
accurate model is
\begin{equation}
y \,=\, A f \,+\,e~, \label{eq:2}
\end{equation}
in which the data are corrupted by an (unknown) noise. In order to
deal with this type of reconstruction problem a {\it regularization}
mechanism is required \cite{enhane96}. Regularization techniques try,
as much as possible, to take advantage of (often vague) prior
knowledge one may have about the nature of $f$. The approach in this
paper is tailored to the case when $f$ can be represented by a
{\it sparse} expansion, i.e.,
when $f$ can be represented by a series
expansion with respect to an orthonormal basis or a
frame \cite{D2,ch03-1} that has only a small number of large
coefficients. In this paper, as in \cite{dadede04}, we model the sparsity
constraint by adding an $\ell_1-$term to a functional to be
minimized; it was shown in \cite{dadede04} that
this assumption does indeed correspond to a regularization scheme.


%

Several types of signals appearing in nature admit sparse frame
expansions and thus, sparsity is a realistic assumption for a very
large class of problems. For instance, natural images are well
approximated by sparse expansions with respect to wavelets or
curvelets \cite{D2,cado04}.

Sparsity has had already a long history of successes.
The design of frames 
for sparse representations of digital signals  has led to extremely
efficient compression methods,
such as JPEG2000 and MP3 
\cite{ma99}.
A new generation of optimal numerical
schemes has been developed for the computation of sparse solutions of
differential and integral equations, exploiting adaptive and greedy
strategies \cite{coh03,CDD1,CDD2,DFR,DFRSW}. The use of sparsity in
inverse problems for data recovery is the most recent step of this concept's
long career of ``simplifying and understanding complexity'', with an
enormous potential in applications \cite{an05,caur04,
cohore04,dama03, dadede04, date05, do92,do95,do95-1,fopi07,
fora06-1,fora07,LoNoDaDa, rate05,te05}. In particular, the
observation that it is possible to reconstruct sparse signals from
vastly incomplete information just seeking for the $\ell_1$-minimal
solutions \cite{cataXX,carotaXX,do04,ra05-7} has led to a new
line of research called {\it sparse recovery} or {\it compressed
sensing}, with very fruitful mathematical and applied results.

\section{Framework and Notations}

Before starting our discussion let us briefly introduce some of the
notations we will need. For some countable index set $\Lambda$
we denote by $\ell_p=\ell_p(\Lambda)$, $1 \leq p \leq \infty$, the
space of real sequences $x=(x_\lambda)_{\lambda \in \Lambda}$ with
norm
\begin{displaymath}
\|x\|_p := \left(\sum_{\lambda \in \Lambda}
|x_\lambda|^p\right)^{1/p}, \quad 1\leq p < \infty
\end{displaymath}
and $\|x\|_\infty \,:=\, \sup_{\lambda \in \Lambda} |x_\lambda|$ as
usual. For simplicity of notation, in the following $\| \cdot\|$
will denote the $\ell_2$-norm $\|\cdot\|_2$.\\
As 
is customary for an index set, we assume we have a natural
enumeration order for the elements of $\Lambda$, using (implicitly)
a one-to-one map $\mathcal{N}$ from $\Lambda$ to $\N$. In some convergence
proofs, we shall use the shorthand notations $|\lambda|$ for $\mathcal{N}(\lambda)$,
and (in the case where $\Lambda$ is infinite) $\lambda \rightarrow \infty$
for $\mathcal{N}(\lambda) \rightarrow \infty$.\\
We also assume that we have a suitable frame $\{\psi_\lambda:
\lambda \in \Lambda\} \subset \mathcal K$ indexed by the countable
set $\Lambda$. This means that there exist constants $c_1,c_2 > 0$
such that
\begin{equation}\label{frame_ineq}
c_1 \|f\|^2_\mathcal K \leq \sum_{\lambda \in \Lambda} |\langle f,
\psi_\lambda\rangle |^2 \leq c_2 \|f\|_\mathcal K^2,\qquad \mbox{ for
all } f \in \mathcal K.
\end{equation}
Orthonormal bases are particular examples of frames, but there also
exist many interesting frames in which the $\psi_{\lambda}$ are not
linearly independent.  Frames allow for a (stable) series expansion
of any $f \in \mathcal K$ of the form
\begin{equation}\label{expand_f}
f \,=\, \sum_{\lambda \in \Lambda} x_\lambda
\psi_\lambda \,=:\,\mbox{F} x \,,
\end{equation}
where $x = (x_\lambda)_{\lambda \in \Lambda} \in \ell_2(\Lambda)$.
The linear operator $\mbox{F} : \ell_2(\Lambda) \to \mathcal K$
(called the {\it synthesis map} in frame theory) is bounded
because of
(\ref{frame_ineq}). When
$\{\psi_\lambda: \lambda \in \Lambda\}$ is a frame but not a basis,
the coefficients $x_\lambda$ need not be unique. For more details on
frames and their differences from bases we refer to \cite{ch03-1}.

We shall assume that  $f$ is sparse,  i.e., that $f$ can be written
by a series of the form (\ref{expand_f}) with only a small number of
non-vanishing coefficients $x_\lambda$ with respect to the frame
$\{\psi_\lambda\}$, or that $f$ is {\em compressible}, i.e., that $f$
can be well-approximated by such a sparse expansion. This can be
modeled by assuming that the sequence $x$ is contained in a
(weighted) $\ell_1(\Lambda)$-space. Indeed, the minimization of the
$\ell_1(\Lambda)$ norm promotes such sparsity. (This has been known
for many years, and put to use in a wide range of applications, most
notably in statistics. David Donoho calls one form of it the Logan phenomenon in
\cite{DoSt89} -- see also \cite{DoLo92} --, after its first observation
by Ben Logan \cite{Lo65}.)
These considerations lead us to model the reconstruction of a sparse
$f$ as the minimization of the following functional:
\begin{equation}
\displaystyle F_\tau(x)=\|Kx-y\|^2_{\mathcal H}+2\tau
\|x\|_1,\label{functional}
\end{equation}
where we will assume that the data $y$ and the linear operator $K:=
A \circ \mbox{F}:\ell_2(\Lambda) \rightarrow \mathcal{H}$ are given. The
second term in (\ref{functional}) is often called the {\em
penalization} or {\em regularizing} term; the first term goes by the
name of {\em discrepancy},
\begin{equation}
D(x):=\|Kx-y\|^2_{\mathcal H}.
\label{discrepancy}
\end{equation}
In what follows we shall drop the subscript ${\mathcal H}$, because
the space in which we work will always be clear from the context. We
discuss the problem of finding (approximations to) $\bar x(\tau)$ in
$\ell_2(\Lambda)$ that minimize the functional (\ref{functional}).
(We adopt the usual convention that for $u \in \ell_2(\Lambda)
\setminus \ell_1(\Lambda)$, the penalty term ``equals'' $\infty$,
and that, for such $u$, $F_\tau(u) > F_\tau(x)$ for all $x\in
\ell_1(\Lambda)$. Since we want to minimize $F_\tau$, we shall
consider, implicitly, only $x \in \ell_1(\Lambda)$.)  The solutions
$\bar f(\tau)$ to the original problem are then given by $\bar
f(\tau) = \mbox{F} \bar x(\tau)$.

Several authors have proposed independently an iterative
soft-threshold\-ing algorithm to approximate the solution $\bar
x(\tau)$ \cite{fino03,stcado03,stngmu03,efhajoti04}. More precisely,
$\bar x(\tau)$ is the limit of sequences $x^{(n)}$ defined
recursively by
\begin{equation}
x\np1=\S_\tau\left[x\n+K^* y-K^* K x\n\right]~~,
\label{DDD04iteration}
\end{equation}
starting from an arbitrary $x^{(0)}$, where $\S_\tau$ is the
soft-thresholding operation defined by
$\S_\tau(x)_\lambda=S_\tau(x_\lambda)$ with
\begin{equation}
S_\tau(x)=\left\{\begin{array}{lll}x-\tau & & x>\tau\\
0 & & |x|\leq \tau\\
x+\tau & & x< -\tau\end{array}\right. .
\label{soft-thresholding}
\end{equation}
Convergence of this algorithm was proved in \cite{dadede04}.
Soft-thresholding plays a role in this problem because it leads to
the unique minimizer of a functional combining $\ell_2$ and
$\ell_1-$norms, i.e., (see \cite{ChDeLeLu,dadede04})
\begin{equation}
\S_\tau(a)=\arg\!\min_{x \in \ell_2(\Lambda)}\left( \|x-a\|^2+2\tau
\|x\|_1\right).
\end{equation}
We will call the iteration \eqref{DDD04iteration} the
\emph{iterative soft-thresholding algorithm} or the
\emph{thresholded Landweber iteration}.

\begin{figure}[h]
\begin{center}
\resizebox{\textwidth}{!}{\includegraphics{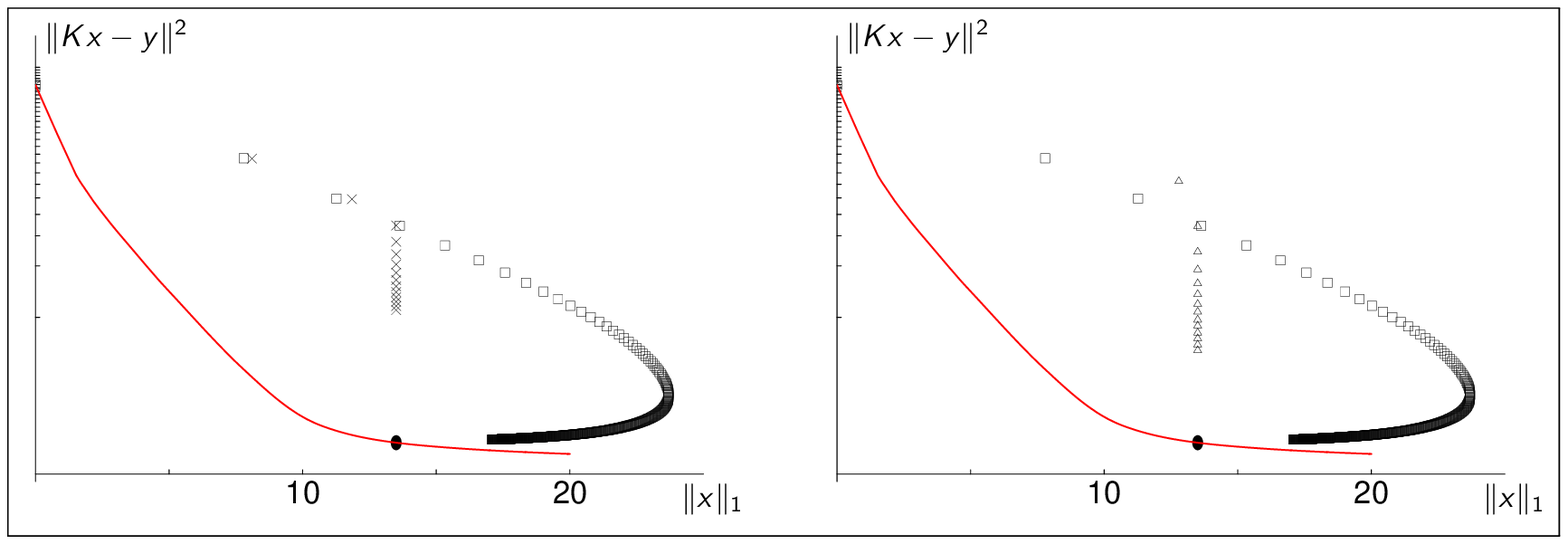}}
\hspace{1 cm} (a) \hspace{8 cm} (b)
\end{center}
\caption{The path, in the $\|x\|_1$ vs. $\|Kx-y\|^2$ plane, followed
by the iterates $x\n$ of three different iterative algorithms. The
operator $K$ and the data $y$ are taken from a seismic tomography
problem \cite{LoNoDaDa} (see also Section \ref{numericsection}). The
boxes (in both (a) and (b))
correspond to the thresholded Landweber algorithm. In this example,
iterative thresholded Landweber
(\ref{DDD04iteration}) first overshoots the $\ell_1$ norm of the
limit (represented by the fat dot), and then requires a large number of
iterations to reduce $\|x\n\|_1$ again (500 are shown
in this figure). In (a)
the crosses correspond to the path followed by the
iterates of the projected Landweber iteration (\ref{PLiteration});
in (b) the triangles correspond to the projected
steepest descent iteration (\ref{PGiteration}); in both cases,
only 15 iterates are shown.
The discrepancy decreases more quickly for projected steepest descent
than for the projected
Landweber algorithm. How this translates into faster convergence (in
norm) is discussed in Section \ref{numericsection}. The solid line
corresponds to the limit {\it trade-off curve}, generated by $\bar
x(\tau)$ for decreasing values of $\tau>0$. The vertical axes uses a
logarithmic scale for clarity.}\label{pathsfig1}
\end{figure}

\section{Discussion of the Thresholded Landweber Iteration}

The problem of finding the sparsest solution to the under-determined
linear equation $Kx=y$ is a hard combinatorial problem, not
tractable numerically except in relatively low dimensions.  For some
classes of $K$, however, one can prove that the problem reduces to
the convex optimization problem of finding the solution with the
smallest $\ell_1$ norm \cite{do04,cataXX,carota06-1,carotaXX}. Even
for $K$ outside this class, $\ell_1-$ minimization seems to lead to
very good approximations to the sparsest solutions. It is in this
sense that an algorithm of type (\ref{DDD04iteration}) could
conceivably be called `fast': it is fast compared to a brute-force
exhaustive search for the sparsest $x$.

A more honest evaluation of the speed of convergence of algorithm
(\ref{DDD04iteration}) is a comparison with
\emph{linear} solvers that minimize the corresponding $\ell_2$
penalized functional, such as, e.g., the conjugate gradient method.
One finds, in practice, that the thresholded Landweber iteration
(\ref{DDD04iteration}) is not competitive at all in this comparison.
It is, after all, the composition of thresholding with the (linear)
Landweber iteration $x\np1=x\n+K^* y-K^* K x\n$, which is a gradient
descent algorithm with a fixed step size, known to converge usually
quite slowly; interleaving it with the nonlinear thresholding
operation does unfortunately not change this slow convergence. On
the other hand, this nonlinearity did foil our attempts to ``borrow
a leaf'' from standard linear steepest descent methods by using an
adaptive step length -- once we start taking larger steps, the
algorithm seems to no longer converge in at least some numerical
experiments.

We take a closer look at the characteristic dynamics of the
thresholded Landweber iteration in Figure \ref{pathsfig1}. As this
plot of the discrepancy
$\mathcal{D}(x\n)=\|Kx^{(n)}-y\|^2$
versus $\|x^{(n)}\|_1$ shows, the algorithm converges initially
relatively fast, then it overshoots the value $\|\bar{x}(\tau)\|_1$
(where $\bar{x}(\tau):=\lim_{n \rightarrow \infty}x^{(n)}$), and it
takes very long to re-correct back. In other words, starting from
$x^{(0)}=0$, the algorithm generates a path $\{x^{(n)}; \,n \in
\mathbb{N}\}$ that is initially fully contained in the $\ell_1$-ball
$B_R:=\{x \in \ell_2(\Lambda); \|x\|_1 \leq R\}$, with $R:=\|\bar
x(\tau)\|_1$. Then it gets out of the ball to slowly inch back to it
in the limit. A first intuitive way to avoid this long ``external''
detour is to force the successive iterates to remain within the ball
$B_R$. One method to achieve this is to substitute for the
thresholding operations the  projection $\P_{B_R}$, where, for any
closed convex set $C$, and any $x$, we define $\P_C(x)$ to be the unique
point in $C$ for which the $\ell_2-$distance to $x$ is minimal. With
a slight abuse of notation, we shall denote $\P_{B_R}$ by $\P_R$;
this will not cause confusion, because it will be clear from the
context whether the subscript of $\P$ is a set or a positive number.
We thus obtain the following algorithm: Pick an arbitrary $x^{(0)} \in \ell_2(\Lambda)$, for example $x^{(0)}=0$, and iterate
\begin{equation}
x\np1=\P_R \left[x\n+K^* y-K^* K x\n\right]. \label{PLiteration}
\end{equation}
We will call this the {\it projected Landweber iteration}.

The typical dynamics of this projected
Landweber algorithm are illustrated in Fig. \ref{pathsfig1}(a).
The norm $\|x^{(n)}\|_1$ no longer overshoots
$R$, but quickly takes on the limit value (i.e., $\|\bar
x(\tau)\|_1$); the speed of convergence remains very slow, however.
In this projected Landweber iteration case, modifying the iterations
by introducing an adaptive ``descent parameter'' $\beta^{(n)}>0$ in
each iteration, defining $x^{(n+1)}$ by
\begin{equation}
x\np1=\P_R \left[x\n+\beta^{(n)}K^* (y-K x\n) \right],
\label{PGiteration}
\end{equation}
does lead, in numerical simulations, to promising, converging
results (in which it differs from the soft-thresholded Landweber iteration, where introducing such a descent parameter did not lead to numerical convergence, as noted above).

The typical dynamics of this modified algorithm are illustrated in Fig. \ref{pathsfig1}(b),
which clearly shows the larger steps and faster convergence (when compared
with the projected Landweber iteration in Fig. \ref{pathsfig1}(a)).
We shall refer to this modified algorithm as the {\it
projected gradient iteration} or the {\it projected steepest
descent}; it will be the main topic of this paper.

The main issue is to determine how large we can choose the
successive $\beta^{(n)}$, and still prove norm convergence of the
algorithm in $\ell_2(\Lambda)$.

There exist results in the literature on convergence of projected
gradient iterations, where the projections are (as they are here)
onto convex sets, see, e.g., \cite{AlIsSo98,CW} and references therein.
These results treat iterative projected gradient methods in much
greater generality than we need: they allow more general functionals
than $\mathcal{D}$, and the convex set on which the iterative
procedure projects need not be bounded. On the other hand, these
general results typically have the  following restrictions:
\begin{itemize}
\item The convergence in infinite-dimensional Hilbert spaces (i.e.,
$\Lambda$ is countable but infinite) is proved only in the weak
sense and often only for subsequences;
\item In \cite{AlIsSo98} the descent parameters are typically restricted to cases for which
$\lim_{n \to \infty}  \beta^{(n)} =0$. In \cite{CW}, it is
shown that the algorithm converges weakly for any choice of $\beta^{(n)} \in \left [ \varepsilon, \frac{2-\varepsilon}{\|K\|} \right ]$, for $\varepsilon>0$ arbitrarily small. Of most interest to us is the case where the $\beta^{(n)}$ are picked adaptively, can grow with $n$, and are not limited to values below $\frac{2}{\|K\|}$; this case
is  not covered by the methods of either \cite{AlIsSo98} or \cite{CW}.
\end{itemize}

To our knowledge there are no results in the literature for which
the whole sequence $(x^{(n)})_{n \in \mathbb{N}}$ converges in the
Hilbert space norm to a unique accumulation point, for ``descent
parameters'' $\beta^{(n)} \geq 2$.
It is worthwhile emphasizing that strong convergence is not
automatic: in \cite[Remark 5.12]{CW}, the authors provide a counterexample in
which strong convergence fails. (This question had been open for some time.)
One of the main results of this
paper is to prove a theorem that establishes exactly this type of
convergence; see Theorem \ref{finaltheorem} below. Moreover, the result is achieved by imposing a choice of $\beta^{(n)} \geq 1$ which ensures a monotone decay of a suitable energy. This establishes a principle of {\it best descent} similar to the well-known steepest-descent in unconstrained minimization.

Before we get to this theorem, we need to build some more machinery
first.

\section{Projections onto $\ell_1$-Balls via Thresholding Operators}

\label{l1proj}

In this section we discuss some properties of $\ell_2$-projections
onto $\ell_1$-balls. In particular, we investigate their relations
with thresholding operators and their explicit computation. We also
estimate the time complexity of such projections in finite
dimensions.

We first observe a useful property of the soft-thresholding
operator.

\begin{Lemma}
For any fixed $a\in \ell_2(\Lambda)$ and for $\tau>0$,
$\|\S_\tau(a)\|_1$ is a piecewise linear, continuous, decreasing
function of $\tau$; moreover, if $a \in \ell_1(\Lambda)$ then $\|\S_{0}(a)\|_1
=\|a\|_1$ and $\|\S_{\tau}(a)\|_1=0$
for $\tau\geq\max_i|a_i|$. \label{piecewiselemma}
\end{Lemma}
{\em Proof:} $\|\S_\tau(a)\|_1=\sum_\lambda|S_\tau(a_\lambda)|=\sum_\lambda
S_\tau(|a_\lambda|)=\sum_{|a_\lambda|>\tau}(|a_\lambda|-\tau)$;  the sum in the right
hand side is finite for $\tau>0$. 
\hfill $\Box$

\begin{figure}
\begin{center}
\includegraphics{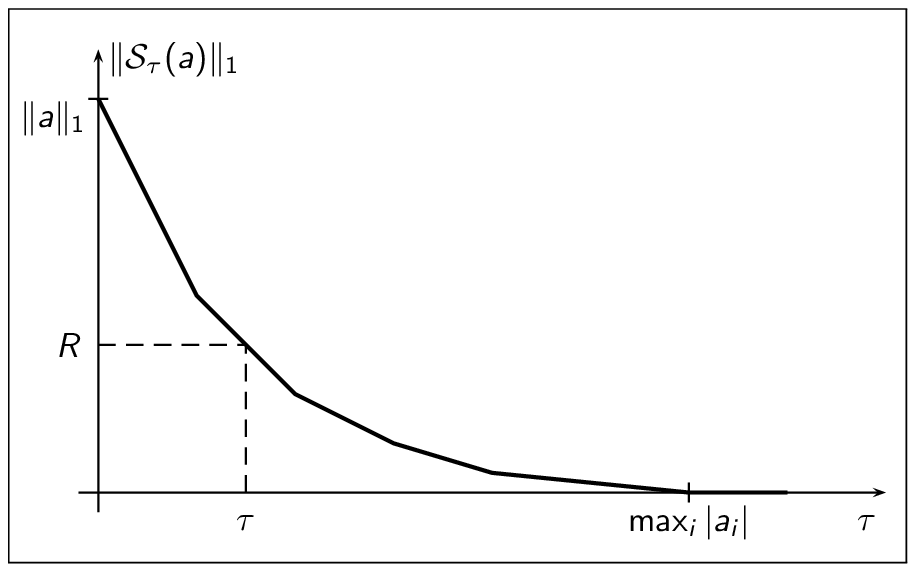}
\end{center}
\caption{For a given vector $a\in\ell_2$, $\|\S_\tau(a)\|_1$ is a
piecewise linear continuous and decreasing function of $\tau$
(strictly decreasing for
$\tau<\max_i|a_i|$) . The knots are located at
$\{|a_i|, i:1\ldots m\}$ and $0$. Finding $\tau$ such that
$\|\S_\tau(a)\|_1=R$ ultimately comes down to a linear
interpolation. The figure is made for the finite dimensional
case.}\label{normfig}
\end{figure}

A schematic illustration is given in Figure \ref{normfig}.

The following lemma shows that the $\ell_2$ projection $\P_R(a)$ can
be obtained by a suitable thresholding of $a$.

\begin{Lemma}
If $\|a\|_1>R$, then the $\ell_2$ projection of $a$ on the $\ell_1$
ball with radius $R$ is given by $\P_R(a)=\S_\mu(a)$ where $\mu$
(depending on $a$ and $R$) is chosen such that $\|\S_\mu(a)\|_1=R$.
If $\|a\|_1\leq R$ then $\P_R(a)= \S_0(a)=
a$.\label{l1projectionlemma}
\end{Lemma}
{\em Proof:}  Suppose $\|a\|_1>R$. Because, by Lemma
\ref{piecewiselemma}, $\|\S_\mu(a)\|_1$ is continuous in $\mu$  and
$\|\S_\mu(a)\|_1=0$ for sufficiently large $\mu$, we can choose
$\mu$ such that $\|\S_\mu(a)\|_1=R$. (See Figure \ref{normfig}.) On
the other hand (see above, or \cite{ChDeLeLu,dadede04}), $b=\S_\mu(a)$ is the
unique minimizer of $\|x-a\|^2+2\mu\|x\|_1$, i.e.,
\begin{displaymath}
\|b-a\|^2+2\mu\|b\|_1< \|x-a\|^2+2\mu\|x\|_1
\end{displaymath}
for all $x\neq b$.  Since $\|b\|_1=R$, it follows that
\begin{displaymath}
\forall x\in B_R,\,\, x \neq b \,:\qquad\qquad  \|b-a\|^2< \|x-a\|^2
\end{displaymath}
Hence $b$ is closer to $a$ than any other $x$ in $B_R$. In other
words, $\P_R(a)=b=\S_\mu(a)$. \hfill $\Box$

These two lemmas prescribe the following simple recipe for computing
the projection $\P_R(a)$. In a first step, sort the absolute values
of the components of $a$ (an $\mathcal{O}(m\log m)$ operation if $\#
\Lambda = m$ is finite), resulting in the rearranged sequence
$\left(a^{\ast}_{\ell}\right)_{\ell=1,\dots,m}$, with
$a^{\ast}_{\ell} \geq a^{\ast}_{\ell+1 }\geq 0$ for all $\ell$.
Next, perform a search to find $k$ such that
\begin{displaymath}
\|\S_{a^{\ast}_k}(a)\|_1=\sum_{\ell=1}^{k-1}
\left(a^{\ast}_{\ell}-a^{\ast}_{k}\right) \leq R < \sum_{\ell=1}^{k}
\left(a^{\ast}_{\ell}-a^{\ast}_{k+1}\right)=
\|\S_{a^{\ast}_{k+1}}(a)\|_1
\end{displaymath}
or equivalently,
\begin{displaymath}
\|\S_{a^{\ast}_k}(a)\|_1=\sum_{\ell=1}^{k-1}\,\ell\,
\left(a^{\ast}_{\ell}-a^{\ast}_{\ell+1}\right) \leq R
<\sum_{\ell=1}^{k}\ell\left(a^{\ast}_{\ell}-a^{\ast}_{\ell+1}\right)=
\|\S_{a^{\ast}_{k+1}}(a)\|_1;
\end{displaymath}
the complexity of this step is again $\mathcal{O}(m\log m)$.
Finally, set \\$\nu:=k^{-1}\left(R-\|\S_{a^{\ast}_k}(a)\|_1\right)$,
and $\mu:=a^{\ast}_k+\nu$. Then
\begin{eqnarray*}
\|\S_{\mu}(a)\|_1&=& \sum_{i \in \Lambda} \max(|a_{i}|-\mu,0)
=  \sum_{\ell=1}^k \left(a^{\ast}_{\ell} - \mu\right) \\
&=&  \sum_{\ell=1}^{k-1} \left(a^{\ast}_{\ell} - a^{\ast}_{k}
\right) \,+\, k \nu \,=\, \|\S_{a^{\ast}_k}(a)\|_1 + k \nu = R.
\end{eqnarray*}

\noindent These formulas were also derived in \cite[Lemma 4.1 and
Lemma 4.2]{fora06-1}, by observing that $\P_R(a) = a - \mathbb S_{
R}^\infty(a)$, where
\begin{equation}\label{def_Svq}
\mathbb S^{\infty}_{R}(a) \,=\, \arg\!\min_{x\in \mathbb R^m}
(\|x-a\|^2 +  2 R \|x\|_\infty),\quad x\in \mathbb R^m.
\end{equation}
The latter is again a thresholding operator, but it is related to an
$\ell_\infty$ penalty term. Similar descriptions of the $\ell_2$
projection onto $\ell_1$ balls appear also in \cite{CaRo04}.

Finally, $\P_R$ has the following additional properties:
\begin{Lemma}\label{convex_basic_lm}
For any $x \in \ell_2(\Lambda)$, $\P_R(x)$ is characterized as the unique vector in
$B_R$ such that
\begin{equation}
\langle w - \P_R(x) , x - \P_R(x) \rangle \leq 0, \text { for all }
w \in B_R.\label{convex_basic}
\end{equation}
Moreover the projection $\P_R$ is non-expansive:
\begin{equation}
\|\P_R(x)-\P_R(x')\|\leq\|x-x'\|\label{contraction}
\end{equation}
for all $x,y\in \ell_2(\Lambda)$.\label{contractionlemma}
\end{Lemma}

\noindent The proof is standard for projection operators onto convex
sets; we include it because its technique will be used often in this
paper.

\noindent {\em Proof:} Because $B_R$ is convex, $(1-t)\P_R(x)+t\,
w\in B_R$ for all $w \in B_R$ and $t \in [0, 1]$. It follows that
$\|x-\P_R(x)\|^2\leq \|x-[(1-t)\,\P_R(x)+t \,w] \,\|^2$ for all $t
\in [0, 1]$. This implies
\begin{displaymath}
0\leq -2 t\, \langle w-\P_R(x),x-\P_R(x)\rangle+t^2 \, \|w-\P_R(x)\|^2
\end{displaymath}
for all $t\in [0,1]$. It follows that
\begin{displaymath}
\langle w-\P_R(x),x-\P_R(x)\rangle\leq 0 ~,
\end{displaymath}
which proves (\ref{convex_basic}).

\noindent Setting $w\,=\,\P_R(x')$ in (\ref{convex_basic}), we get,
for all $x, x'$,
\begin{displaymath}
\langle \P_R(x')-\P_R(x),x-\P_R(x)\rangle\leq 0
\end{displaymath}
Switching the role of $x$ and $x'$ one finds:
\begin{displaymath}
\langle \P_R(x')-\P_R(x),x'-\P_R(x')\rangle\geq 0
\end{displaymath}
By combining these last two inequalities, one finds:
\begin{displaymath}
\langle \P_R(x')-\P_R(x),x'-x-\P_R(x')+\P_R(x)\rangle\geq 0
\end{displaymath}
or
\begin{displaymath}
\|\P_R(x')-\P_R(x)\|^2\leq \langle \P_R(x')-\P_R(x),x'-x\rangle  \,
;
\end{displaymath}
by Cauchy-Schwarz this gives
\begin{displaymath}
\|\P_R(x')-\P_R(x)\|^2\leq \langle \P_R(x')-\P_R(x),x'-x\rangle \leq
\| \P_R(x')-\P_R(x)\| \|x'-x\| \, ,
\end{displaymath}
from which inequality (\ref{contraction}) follows.\hfill $\Box$

\section{The Projected Gradient Method}

\label{proofsection}

We have now collected all the terminology needed to identify some
conditions on the $\beta^{(n)}$ that will ensure convergence of the
$x^{(n)}$, defined by \eqref{PGiteration}, to $\tilde{x}_R$, the
minimizer in $B_R$ of $\mathcal{D}(x)=\|Kx-y\|^2$. For
notational simplicity we set $r\n=K^*(y- K x\n)$. With this
notation, the thresholded Landweber iteration (\ref{DDD04iteration})
can be written as
\begin{equation}
x\np1=\S_\tau\left(x\n+r\n\right).
\label{threshLWiteration}
\end{equation}
As explained above, we consider, instead of straightforward
soft-thresholding with fixed $\tau$,  adapted soft-thresholding
operations $\S_{\mu(R,x^{(n)}+r^{(n)})}$ that correspond to the
projection operator $\P_R$:
\begin{equation}
x\np1=\P_R\left(x\n+r\n\right).\label{projLWiteration}
\end{equation}
The dependence of $\mu(R,x^{(n)}+r^{(n)})$ on $R$ is described
above; $R$ is kept fixed throughout the iterations. If, for a given
value of $\tau$, $R$ were picked such that $R=
R_{\tau}:=\|\bar{x}_{\tau}\|_1$ (where $\bar{x}_{\tau}$ is the
minimizer of $\|Kx-y\|^2+2 \tau \|x\|_1$), then Lemma
\ref{l1projectionlemma} would ensure that $\bar{x}_{\tau}=
\tilde{x}_R$. Of course, we don't know, in general, the exact value
of $\|\bar{x}_{\tau}\|_1$, so that we can't use it as a guideline to
pick $R$. In practice, however, it is customary to determine
$\bar{x}_{\tau}$ for a range of $\tau$-values; this then amounts to
the same as determining $\tilde{x}_R$ for a range of $R$.

\noindent We now propose to change the step $r\n$ into a step
$\beta\n r\n$ (in the spirit of the ``classical'' steepest descent
method), and to define the algorithm: Pick an arbitrary $x^{(0)} \in \ell_2(\Lambda)$, for example $x^{(0)}=0$, and iterate
\begin{equation}
x\np1=\P_R\left(x\n+\beta\n r\n\right).\label{projSDiteration}
\end{equation}
In this section we prove the norm convergence of this algorithm to a
minimizer $\tilde{x}_R$ of $\|Kx-y\|^2$ in $B_R$,
under some assumptions on the descent parameters $\beta\n \geq 1$.

\subsection{General properties}

\label{genpropsec}

\noindent We begin with the following characterization of the
minimizers of $\mathcal{D}$ on $B_R$.
\begin{Lemma}\label{fixpt}
The vector $\tilde{x}_R \in \ell_2(\Lambda)$ is a minimizer of
$\mathcal{D}(x)=\| K x -y\|^2$ on $B_R$ if and only if
\begin{equation} \label{fixpteq}
\P_R(\tilde{x}_R + \beta K^*(y - K \tilde{x}_R)) = \tilde{x}_R,
\end{equation}
for any $\beta >0$, which in turn is equivalent to the requirement
that
\begin{equation}
\langle K^*(y - K \tilde{x}_R), w - \tilde{x}_R \rangle \leq 0,
\text { for all } w \in B_R.
\end{equation}
\end{Lemma}
To lighten notation, we shall drop the subscript $R$ on
$\tilde{x}_R$
whenever no confusion is possible.\\
{\em Proof:} If $\tilde x$ minimizes $\mathcal{D}$ on $B_R$, then
for all $w \in B_R$, and for all $t \in [0,1]$,
\begin{eqnarray*}
\mathcal{D}(\tilde x) &\leq& \mathcal{D}( (1-t) \tilde x + t w),
\text { or } \\
\| K \tilde x - y\|^2 &\leq& \| K \tilde x - y +
t K( w - \tilde x)\|^2 \text{, or }\\
 0&\leq&  2 t \langle K  \tilde x - y, K(w -\tilde x)
\rangle + t^2 \|  K(w - \tilde x)\|^2.
\end{eqnarray*}
This implies
\begin{equation}
 \langle  K^*(y -K \tilde x), w -\tilde x \rangle \leq 0.
\end{equation}
It follows from this that, for all $w \in B_R$ and for all $\beta
>0$,
\begin{equation}
 \langle  \tilde x + \beta K^*(y -K \tilde x) - \tilde x , w -\tilde x \rangle \leq 0,
\end{equation}
By Lemma \ref{contractionlemma} this implies \eqref{fixpteq}.


Conversely, if $\P_R(\tilde x + \beta K^*(y - K \tilde x)) = \tilde
x$, then for all $w \in B_R$ and for all $t \in [0,1]$:
\begin{eqnarray*}
\| (\tilde x + \beta  K^*(y - K \tilde x)) - ((1 -t) \tilde x +t
w)\|^2
&\geq& \| (\tilde x + \beta  K^*(y - K \tilde x)) - \tilde x\|^2,\\
\text{ or } \| \beta  K^*(y - K \tilde x) + t(\tilde x -w)\|^2
&\geq& \| \beta  K^*(y - K \tilde x))\|^2,\\
\Rightarrow \quad 2 t \beta \langle  K^*(y - K \tilde x), \tilde x
-w \rangle + t^2 \| \tilde x -w\|^2 &\geq& 0.
\end{eqnarray*}
This implies
\begin{displaymath}
\langle  K^*(y - K \tilde x), \tilde x -w \rangle \geq 0 \qquad
\mathrm{or}\qquad \langle  y - K \tilde x, K(\tilde x -w) \rangle
\geq 0.
\end{displaymath}
In other words:
\begin{displaymath}
- \| y - K \tilde x\|^2 - \| K \tilde x - K
w\|^2 + \| (y - K \tilde x) + K(\tilde x -
w)\|^2\geq 0
\end{displaymath}
or
\begin{displaymath}
\mathcal{D}(\tilde x) + \|K(\tilde x -w)\|^2 \leq
\mathcal{D}(w).
\end{displaymath}
This implies that $\tilde x$ minimizes $\mathcal{D}$ on $B_R$.
\hfill $\Box$\\

The minimizer of $\mathcal{D}$ on $B_R$ need not be unique. We have,
however
\begin{Lemma}\label{kern}
If $\tilde x, \tilde{\tilde{x}}$ are two distinct minimizers of
$\mathcal{D}(x)=\|Kx-y\|^2$ on $B_R$, then $K \tilde x = K
\tilde{\tilde{x}}$, i.e., $\tilde x - \tilde{\tilde{x}} \in \ker K$.
\\ Conversely, if $\tilde x, \tilde{\tilde{x}} \in B_R$, if $\tilde
x$ minimizes $\|Kx-y\|^2$ and if $\tilde x - \tilde{\tilde{x}} \in
\ker K$ then $\tilde{\tilde{x}}$ minimizes $\|Kx-y\|^2$ as well.
\end{Lemma}
{\em Proof:} The converse is obvious; we prove only the direct
statement.
From the last inequality in the proof of Lemma \ref{fixpt} we
obtain
$\mathcal{D}(\tilde x) + \|K(\tilde x -\tilde{\tilde x})\|^2 \leq
\mathcal{D}(\tilde{\tilde x}) = \mathcal{D}(\tilde x)
$,
which implies $\|K(\tilde{\tilde x} -{\tilde x})\|=0$.
\hfill $\Box$\\

In what follows we shall assume that the minimizers of $\mathcal{D}$
in $B_R$ are not global minimizers for $\mathcal{D}$, i.e., that
$K^*(y- K \tilde x) \neq 0$. We know from Lemma
\ref{l1projectionlemma} that $\mathbb P_R(a)$ can be computed for
$\|a\|_1 > R$ simply by finding the value $\mu > 0$ such that $\|
\mathbb S_\mu(a)\|_1=R$; one has then $\P_R (a) =    \mathbb
S_\mu(a)$. Using this we prove

\begin{Lemma}
\label{thrval}
Let $u$ be the common image under $K$ of all minimizers of
$\mathcal{D}$ on $B_R$, i.e., for all $\tilde x$ minimizing
$\mathcal{D}$ in $B_R$, $K \tilde x =u$. Then there exists a unique
value $\tau>0$ such that, for all $\beta >0$ and for all minimizing
$\tilde x$
\begin{equation}
\mathbb P_R(\tilde x + \beta K^* (y - u)) = \mathbb S_{\tau \beta} (
\tilde x + \beta K^* (y- u)).
\end{equation}
Moreover, for all $\lambda \in \Lambda$ we have that
if there exists a minimizer $\tilde x$ such that $\tilde x_\lambda
\neq 0$, then
\begin{equation}
| (K^* (y-u))_\lambda| = \tau.
\end{equation}
\end{Lemma}
{\em Proof:} From Lemma \ref{l1projectionlemma} and Lemma
\ref{fixpt}, we know that for each minimizing $\tilde x$, and each
$\beta >0$, there exists a unique $\mu(\tilde x,\beta)$ such that
\begin{equation}
\tilde x = \mathbb P_R(\tilde x + \beta K^* (y - u)) = \mathbb
S_{\mu(\tilde x, \beta)} ( \tilde x + \beta K^* (y- u)).
\end{equation}
For $\tilde x_\lambda \neq 0$ we have $\tilde x_\lambda = \tilde
x_\lambda + \beta  (K^* (y - u))_\lambda - \mu(\tilde x, \beta)
\sign \, \tilde x_\lambda$; this implies $ \sign \, \tilde x_\lambda
= \sign \, (\tilde x_\lambda + \beta  (K^* (y - u))_\lambda)$ and also that $ | (K^* (y-u))_\lambda| = \frac{1}{\beta} \mu(\tilde
x,\beta)$. If $\tilde x_\lambda =0$ then $ | (K^*
(y-u))_\lambda| \leq \frac{1}{\beta} \mu(\tilde x, \beta)$. It
follows that $\tau := \mu(\tilde x, \beta)/\beta = \| K^* (y-u)\|_\infty$ does not depend on the choice of $\tilde x$. Moreover, if
there is a minimizer $\tilde x$ for which $\tilde x_\lambda \neq 0$,
then $| (K^* (y-u))_\lambda| = \tau$.
\hfill $\Box$\\

\begin{Lemma}
If, for some $\lambda \in \Lambda$, two minimizers $\tilde x,
\tilde{\tilde{x}}$ satisfy $\tilde x_\lambda \neq 0$ and
$\tilde{\tilde{x}}_\lambda \neq 0$, then $\sign \, \tilde x_\lambda
= \sign \, \tilde{\tilde{x}}_\lambda$.
\end{Lemma}
{\em Proof:} This follows from the arguments in the previous proof;
$\tilde x_\lambda \neq 0$ implies $(K^* (y - u))_\lambda= \tau
\,\sign \, \tilde x_\lambda$. Similarly, $\tilde{\tilde{x}}_\lambda
\neq 0$ implies $(K^* (y - u))_\lambda) = \tau \,\sign \,
\tilde{\tilde{ x}}_\lambda$, so that $\sign \, \tilde x_\lambda =
\sign \, \tilde{\tilde{x}}_\lambda$.
\hfill $\Box$\\

This immediately leads to
\begin{Lemma}\label{gammaset}
For all $\tilde x \in B_R$ that minimize $\mathcal{D}$, there are
only finitely many $\tilde x_\lambda \neq 0$. More precisely,
\begin{equation}
\{\lambda \in \Lambda: \tilde x_\lambda \neq 0\}
\subset \Gamma:= \{ \lambda \in \Lambda: | (K^* (y-u))_\lambda| =\|
(K^* (y-u))\|_\infty\}.
\end{equation}
Moreover, if the vector $e$ is defined by
\begin{equation}
e_\lambda = \left\{ \begin{array}{ll}
0, & \lambda \notin \Gamma\\
\sign ((K^* (y-u))_\lambda),& \lambda  \in \Gamma,
\end{array} \right .
\end{equation}
then $\langle \tilde x, e \rangle = R$ for each minimizer $\tilde x$
of $\mathcal{D}$ in $B_R$.
\end{Lemma}
{\em Proof:} We have already proved the set inclusion. Note that,
since $K^*(y-u) \in \ell_2(\Lambda)$, the set $\Gamma$ is
necessarily a finite set. We also have, for each minimizer $\tilde
x$,
\begin{eqnarray*}
\langle \tilde x, e \rangle &=& \sum_{\lambda \in \Gamma} \tilde x_\lambda e_\lambda \\
&=&
\sum_{\lambda \in \Gamma, \,\tilde x_\lambda \neq 0} \tilde x_\lambda \, \sign ((K^* (y-u))_\lambda) \\
&=&  \sum_{\lambda \in \Gamma,\,\tilde x_\lambda \neq 0} \tilde
x_\lambda \, \sign (\tilde x_\lambda) = \|\tilde x\|_1 =R.
\end{eqnarray*}
\hfill $\Box$

\begin{Remark}\label{conven} By changing, if necessary, signs of the canonical basis vectors, we
can assume, without loss of generality, that $e_\lambda = +1$ for
all $\lambda \in \Gamma$. We shall do so from now on.
\end{Remark}

\subsection{Weak convergence to minimizing accumulation points}

We shall now impose some conditions on the $\beta\n$.
We shall see examples in Section \ref{numericsection} where these
conditions are verified.

\begin{Definition}\label{conditionB}
We say that the sequence $\left(\beta\n\right)_{n \in \N}$ {\em
satisfies Condition (B) with respect to the sequence}
$\left(x\n\right)_{n \in \N}$ if there exists  $n_0$ so that:
\begin{eqnarray*}
\mbox{{(B1)}} && \bar{\beta} := \sup \{ \,\beta\n\,;\, n \in \N
\,\} \, < \, \infty \,\quad\mbox{ and }\quad \inf\{ \,\beta\n\,;\, n
\in \N \,\} \,\geq \,1
\, \\
\mbox{{ (B2)}} && \beta\n \| K(x\np1-x\n) \|^2\leq
r \| x\np1 -x\n \|^2 \, \qquad \forall n \geq n_0.
\end{eqnarray*}
\end{Definition}
We shall often abbreviate this by saying that `the $\beta^{(n)}$
satisfy Condition (B)'. The constant $r$ used in this definition is
$r:=\|K^* K\|_{\ell_2 \to \ell_2}  <1$. (We can always
assume, without loss of generality, that $\|K\|_{\ell_2 \to
\mathcal{H}}<1$; if necessary, this can be achieved
by a suitable rescaling of $K$ and $y$.)

Note that the choice
$\beta^{(n)}=1$ for all $n$, which corresponds to the projected
Landweber iteration, automatically satisfies Condition (B); since we
shall show below that we obtain convergence when the $\beta^{(n)}$
satisfy Condition (B), this will then establish, as a corollary,
convergence of the projected Landweber iteration algorithm \eqref{PLiteration} as well.
We shall be interested in choosing, adaptively, larger values of
$\beta^{(n)}$; in particular, we like to choose $\beta^{(n)}$ as
large as possible.

\begin{Remark}
\begin{itemize}
\item{ Condition (B) is inspired by the standard length-step in the
steepest descent algorithm for the (unconstrained, unpenalized)
functional $\|Kx-y\|^2$. In this case, one can speed up the standard
Landweber iteration $x\np1= x\n + K^*(y -K x\n)$ by defining instead
$x\np1= x\n + \alpha K^*(y -K x\n)$, where $\alpha$ is picked so
that it gives the largest decrease of $\|Kx-y\|^2$ in this
direction. This gives
\begin{equation}
\alpha
 = \left[\|K^*(y -K x\n)\|^2\right]
\left[\|K K^*(y -K x\n)\|^2\right]^{-1}\,.
\end{equation}
In this linear case, one easily checks that $\alpha $ also equals
\begin{equation}
\alpha
 = \left[\|x\np1-x\n\|^2\right]
\left[\|K(x\np1-x\n)\|^2\right]^{-1}\,;
\end{equation}
in fact, it is this latter expression for $\alpha$ (which inspired
the formulation of Condition (B)) that is most useful in proving
convergence of the steepest descent algorithm.}
\item{Because the definition of $x\np1$
involves $\beta^{(n)}$, the inequality (B2), which uses $x\np1$ to impose
a limitation on $\beta^{(n)}$, has an ``implicit'' quality. In practice,
it may not be straightforward to pick $\beta^{(n)}$ appropriately; one could
conceive of trying first a ``greedy'' choice, such as e.g.
$\frac{\|r^{(n)}\|^2}{\|Kr^{(n)}\|^2}$; if this value works, it is retained;
if it doesn't, it can be gradually decreased (by multiplying it with a factor
slightly smaller than 1) until (B2) is satisfied. (A similar way of testing
appropriate step lengths is adopted in \cite{finowr07}.)
}
\end{itemize}
\end{Remark}

In this section we prove that if the sequence $(x\n)_{n \in
\mathbb{N}}$ is defined iteratively by \eqref{PGiteration}, and if
the $\beta^{(n)}$ used in the iteration satisfy Condition (B) (with
respect to the $x\n$), then the (weak) limit of any weakly
convergent subsequence of $(x\n)_{n \in \mathbb{N}}$ is  necessarily
a minimizer of $\mathcal{D}$ in $B_R$.


\begin{Lemma} \label{minimfbeta}
Assume $\|K\|_{\ell_2 \to \mathcal{H}}<1$ and $\beta \geq 1$.
For arbitrary fixed $x$ in $B_R$, define the functional $F_\beta(\cdot;x)$
by
\begin{equation}
\label{funcn}
F_\beta(w;x) := \| K w - y\|^2 - \|K(w -
x)\|^2+ \frac{1}{\beta} \| w-x\|^2\, .
\end{equation}
Then
there is a unique choice for $w$ in $B_R$ that minimizes the
restriction to $B_R$ of $F_\beta(w;x)$. We denote this minimizer by
$T_R(\beta;x)$; it is given by
$T_R(\beta;x)=  \P_R(x + \beta K^*(y -
Kx))$.
\end{Lemma}
{\em Proof:} First of all, observe that
the functional $ F_\beta(\cdot,x)$ is strictly convex, so that it has a
unique minimizer on $B_R$; let $\hat x$ be this  minimizer. Then for
all $w \in B_R$ and for all $t \in [0,1]$
\begin{eqnarray*}
& & \quad  F_\beta( \hat x;x) \leq F_\beta( (1- t)\hat x + t w;x) \\
&\Rightarrow& \quad 2 t \left [ \langle K \hat x - y, K(w -\hat x) \rangle -
\langle K \hat x - K x, K(w - \hat x) \rangle + \frac{1}{\beta} \langle \hat x -
x , w - \hat x \rangle \right ]\\
&&\phantom{xx} + \frac{t^2}{\beta} \| w - \hat x\|^2 \geq 0\\
  &\Rightarrow& \quad  \left [ \beta \langle K x - y, K(w -\hat x) \rangle  +
  \langle \hat x - x , w - \hat x \rangle \right ] + \frac{t }{2} \| w - \hat x\|^2 \geq 0\\
  &\Rightarrow& \quad \langle  \hat x  -x  +\beta K^*(K x -y),w- \hat x  \rangle \geq 0 \\
 &\Rightarrow& \quad \langle x + \beta K^*(y-K x)- \hat x, w- \hat x  \rangle \leq 0.
\end{eqnarray*}
The latter implication is equivalent to $\hat x = \mathbb P_ R( x + \beta K^*(y -
Kx))$ by Lemma \ref{convex_basic_lm}.
\hfill $\Box$\\

\noindent An immediate consequence is

\begin{Lemma}\label{asymptreg}
If the $x\n$ are defined by \eqref{PGiteration}, and the $\beta\n$
satisfy {\em Condition (B)} with respect to the $x\n$, then the
sequence $\left(\mathcal{D}(x\n)\right)_{n \in \mathbb{N}}$ is
decreasing, and
\begin{equation}
\label{asympt} \lim_{n \to \infty} \| x\np1 - x\n \| =0.
\end{equation}
\end{Lemma}
{\em Proof:} Comparing the definition of $x\np1$ in
\eqref{PGiteration} with the statement of Lemma \ref{minimfbeta}, we
see that $x\np1= T_R(\beta\n;x\n)$, so that $x\np1$ is the
minimizer, for $x \in B_R$, of $F_{\beta\n}(x;x\n)$. Setting $\gamma
= \frac{1}{r}-1>0$, we have
\begin{eqnarray*}
\mathcal{D}(x\np1) &\leq & \mathcal{D}(x\np1) + \gamma
\| K(x\np1 - x\n)\|^2 \\
 &= & \| K x\np1 - y\|^2 +(1+\gamma)
 \| K(x\np1 - x\n)\|^2 - \| K(x\np1 - x\n)\|^2 \\
&\leq &  \| K x\np1 - y\|^2 - \| K(x\np1 - x\n)\|^2
+\frac{1}{\beta\n} \|x\np1 - x\n\|^2\\
&=&   F_{\beta\n}(x\np1;x\n) \leq F_{\beta\n}(x\n;x\n) =
\mathcal{D}(x\n).
\end{eqnarray*}
We also have
\begin{eqnarray*}
&&-F_{\beta\np1} ( x\np1;x\np1) + F_{\beta\n} ( x\np1;x\n) \\
&=&
\frac{1}{\beta\n} \| x\np1 - x\n\|^2 - \|K(x\np1 - x\n)\|^2\\
&\geq& \frac{1-r}{\beta\n}  \| x\np1 - x\n\|^2 \geq \frac{1-r}{\bar
\beta}\| x\np1 - x\n\|^2.
\end{eqnarray*}
This implies
\begin{eqnarray*}
\sum_{n=0}^N \| x\np1 - x\n\|^2 &\leq& \frac{\bar \beta}{1-r}
\sum_{n=0}^N \left ( F_{\beta\n} ( x\np1;x\n) -F_{\beta\np1} ( x\np1;x\np1) \right) \\
&\leq& \frac{\bar \beta}{1-r} \sum_{n=0}^N \left ( F_{\beta\n} ( x\n;x\n)
-F_{\beta\np1} ( x\np1;x\np1) \right) \\
&=& \frac{\bar \beta}{1-r} \left ( F_{\beta^{(0)}} ( x^{(0)};x^{(0)})
-F_{\beta^{(N+1)}} ( x^{(N+1)}; x^{(N+1)}) \right) \\
&\leq& \frac{\bar \beta}{1-r} F_{\beta^{(0)}} ( x^{(0)};x^{(0)}).
\end{eqnarray*}
Therefore, the series $\sum_{n=0}^\infty\| x\np1 - x\n\|^2$
converges and $\lim_{n \to \infty} \| x\np1 - x\n \| =0$.
\hfill $\Box$\\

Because the set $\{x\n;n \in \mathbb{N} \}$ is bounded in $\ell_1(\Lambda)$ ($x^n$
are all in $B_R$), it is bounded in
$\ell_2(\Lambda)$ as well (since $\|a\|_2 \leq \|a\|_1$). Because
bounded closed sets in $\ell_2(\Lambda)$ are weakly compact, the
sequence $(x\n)_{n \in \mathbb{N}}$ must have weak accumulation
points. We now have

\begin{Proposition}[Weak convergence to minimizing accumulation points]
\label{weakconv} If $x^\#$ is a weak accumulation point of
$(x\n)_{n \in \mathbb{N}}$ then $x^\#$ minimizes $\mathcal{D}$ in
$B_R$.
\end{Proposition}
{\em Proof:} Let $(x^{(n_j)})_{j \in \N}$ be a subsequence
converging weakly to $x^\#$. Then for all $a \in \ell_2(\Lambda)$
\begin{equation}
\langle K x^{(n_j)} , a \rangle = \langle x^{(n_j)} , K^* a \rangle
\,_{\overrightarrow{ j \to \infty }}\, \langle x^{\#} , K^* a
\rangle = \langle K x^{\#} , a \rangle.
\end{equation}
Therefore $\wlim_{j \to \infty } \, K  x^{(n_j)} = K x^\#$. From
Lemma \ref{asymptreg} we have $\| x\np1 - x\n \|
\,_{\overrightarrow{n \to \infty}}\, 0$, so that we also have
$\wlim_{j \to \infty} \, x^{(n_j+1)} = x^\#$. By the definition of
$x\np1$ ($ x\np1= \P_R ( x\n + \beta\n K^* ( y - K x\n))$), and by
Lemma \ref{contractionlemma}, we have, for all $w \in B_R$,
\begin{equation}
\langle  x\n + \beta\n K^* ( y - K x\n) - x\np1 , w - x\np1 \rangle
\leq 0.
\end{equation}
In particular, specializing to our subsequence and taking the
$\limsup$, we have
\begin{equation}
\limsup_{j \to \infty} \, \langle   x^{(n_j)} -    x^{(n_j+1)} +
\beta^{(n_j)} K^* ( y - K  x^{(n_j)}), w -  x^{(n_j+1)} \rangle \leq
0.
\end{equation}
Because $\| x^{(n_j)} -    x^{(n_j+1)} \| \rightarrow 0$, for $j \to
\infty$, and $ w -  x^{(n_j+1)}$ is uniformly bounded, we have
\begin{equation}
\lim_{j \to \infty} \, | \langle   x^{(n_j)} -    x^{(n_j+1)}, w -
x^{(n_j+1)} \rangle | = 0,
\end{equation}
so that our inequality reduces to
\begin{equation}
\limsup_{j \to \infty} \, \beta^{(n_j)}  \,\langle K^* ( y - K
x^{(n_j)}), w -  x^{(n_j+1)} \rangle \leq 0.
\end{equation}
By adding $\beta^{(n_j)}\langle K^* ( y - K  x^{(n_j+1)}),
x^{(n_j+1)} - x^{(n_j)} \rangle$, which also tends to zero as $j \to
\infty$, we transform this into
\begin{equation}
\limsup_{j \to \infty} \, \beta^{(n_j)}  \,\langle K^* ( y - K
x^{(n_j)}), w -  x^{(n_j)} \rangle \leq 0.
\end{equation}
Since the $\beta^{(n_j)}$ are all in $[1, \bar \beta]$, it follows
that
\begin{equation}
\limsup_{j \to \infty} \, \langle K^* ( y - K  x^{(n_j)}), w -
x^{(n_j)} \rangle \leq 0,
\end{equation}
or
\begin{equation}
\limsup_{j \to \infty} \, \left [ \langle K^* y , w -x^\# \rangle -
\langle K^* K x^\#,w \rangle + \| K x^{(n_j)}\|^2 \right ] \leq 0,
\end{equation}
where we have used the weak convergence of $ x^{(n_j)}$. This can be
rewritten as
\begin{equation}
\label{spec} \langle K^* ( y - K x^\#) , w -x^\# \rangle +
\limsup_{j \to \infty} \, \left [ \| K x^{(n_j)}\|^2 -  \| K
x^{\#}\|^2 \right ] \leq 0.
\end{equation}
Since $\wlim_{j \in \mathbb{N}} \, K  x^{(n_j)} = K x^\#$, we have
$$\limsup_{j \to \infty} \, \left [ \| K x^{(n_j)}\|^2 -
\| K x^{\#}\|^2 \right ] \geq 0.$$ We conclude thus that
\begin{equation}
\langle K^* ( y - K x^\#) , w -x^\# \rangle \leq 0,  \quad \text {
for all } w \in B_R,
\end{equation}
so that $x^\#$ is a minimizer of $\mathcal{D}$ on $B_R$, by Lemma
\ref{fixpt}. \hfill $\Box$

\subsection{Strong convergence to minimizing accumulation points}

In this subsection we show how the weak convergence established in
the preceding subsection can be strengthened into 
norm convergence, again by a series of lemmas. Since the distinction
between weak and strong convergence makes sense only when the index
set $\Lambda$ is infinite, we shall implicitly assume this is the
case throughout this section.

\begin{Lemma}
\label{strongconv} For the subsequence $(x^{(n_j)})_{ j \in
\mathbb{N}}$ defined in the proof of {\em Proposition
\ref{weakconv}},\\ $\lim_{j \to \infty}  K (x^{(n_j)}) = K x^\#$.
\end{Lemma}
{\em Proof:} Specializing the inequality  \eqref{spec} to $w =
x^\#$, we obtain
$$
\limsup_{j \to \infty} \, \left [ \| K x^{(n_j)}\|^2 -  \| K
x^{\#}\|^2 \right ] \leq 0;
$$
together with $\| K x^\#\|^2\leq  \liminf_{j \to \infty} \|K
x^{(n_j)}\|^2$ (a consequence of the weak convergence of
$Kx^{(n_j)}$ to $K x^\#$),
this implies $\lim_{j \to \infty}  \|K
(x^{(n_j)})\|^2 = \|K x^\#\|^2$, and thus $\lim_{j \to \infty}  K
(x^{(n_j)}) = K x^\#$.
\hfill $\Box$\\

\begin{Lemma}
Under the same assumptions as in {\em Proposition \ref{weakconv}},
there exists a subsequence $\left(x^{(n'_{\ell})}\right)_{\ell \in
\N}$ of $(x\n)_{n \in \N}$ such that
\begin{equation}
\lim_{\ell \to \infty} \| x^{(n'_{\ell})}- x^\#\| =0,
\end{equation}
\end{Lemma}
{\em Proof:} Let $(x^{(n_j)})_{ j \in \mathbb{N}}$ be the
subsequence defined in the proof of Proposition \ref{weakconv}.
Define now $u^{(j)} := x^{(n_j)}- x^\#$ and $v^{(j)} := x^{(n_j+1)}-
x^\#$. Since, by Lemma \ref{asymptreg},  $\| x\np1 - x\n \|
\,_{\overrightarrow{n \to \infty}} \,0$,  we have $\| u^{(j)} -
v^{(j)}\| \,_{\overrightarrow{j \to \infty}} \,0$. On the other
hand,
\begin{eqnarray*}
 u^{(j)} - v^{(j)} &=&  u^{(j)} + x^\# - \mathbb P_R \left
 (  u^{(j)} + x^\#  + \beta^{(n_j)} K^* ( y - K(u^{(j)} +x^\#)) \right)\\
&=& u^{(j)} + \mathbb P_R \left ( x^\#  + \beta^{(n_j)} K^* ( y - K x^\#) \right) \\
&& \quad - \mathbb P_R \left (  x^\#  + \beta^{(n_j)} K^* ( y - K
x^\#)+  u^{(j)} -  \beta^{(n_j)} K^* K u^{(j)}\right),
\end{eqnarray*}
where we have used Proposition \ref{weakconv} ($x^\#$ is a
minimizer) and Lemma \ref{fixpt} (so that \\$x^\# = \mathbb P_R \left
( x^\#  + \beta^{(n_j)} K^* ( y - K x^\#) \right)$).  By Lemma
\ref{strongconv}, $\lim_{j \to \infty} \| K u^{(j)} \| =0$. Since
the $\beta^{(n_j)}$ are uniformly bounded, we have, by formula
\eqref{contraction},
\begin{eqnarray*}
&&\!\!\!\!\!\!\!\!\left \| \mathbb P_R \left( x^\#  + \beta^{(n_j)}
K^* ( y - Kx^\#)+  u^{(j)} - \beta^{(n_j)} K^* K u^{(j)}\right)
\right.\\
&&\phantom{xxxxxxxxxxxx} \left . \mathbb P_R \left (  x^\#  + \beta^{(n_j)} K^*
( y - K x^\#)+  u^{(j)}\right) \right \| \\
&& \phantom{xxxxxxxxxxxx} \leq  \beta^{(n_j)} \| K^* K u^{(j)} \|
\,_{\overrightarrow{ j \to \infty }} 0.
\end{eqnarray*}
Combining this with $\| u^{(j)} - v^{(j)}\|\,_{\overrightarrow{ j
\to \infty }} \,0$,  we obtain
$$
\lim_{j \to \infty} \left \| \mathbb P_R \left (  x^\#  +
\beta^{(n_j)} K^* ( y - Kx^\#)+  u^{(j)}\right) \right.
$$
\begin{equation}
\left .
- \mathbb P_R \left
(  x^\#  + \beta^{(n_j)} K^* ( y - Kx^\#)\right) -  u^{(j)} \right
\| =0.
\end{equation}
Since the $ \beta^{(n_j)}$ are uniformly bounded, they must have at
least one accumulation point. Let $\beta^{(\infty)}$ be such an
accumulation point, and choose a subsequence $(j_{\ell})_{\ell \in
\N}$ such that $\lim_{\ell \rightarrow
\infty}\beta^{(n_{j_{\ell}})}= \beta^{(\infty)}$. To simplify
notation, we write $\nell:=n_{j_{\ell}}$, $\uell:= u^{(j_{\ell})}$,
$\vell:= v^{(j_{\ell})}$. We have thus
\begin{equation}\label{helpeq}
\begin{array}{l}
\lim_{\ell \to \infty} \beta^{(
\nell)} = \beta^{(\infty)}\, ,\, \mbox{and } \\
\lim_{\ell \to \infty} \left \| \mathbb P_R \left (  x^\#  +
\beta^{(\nell)} K^* ( y - Kx^\#)+  \uell\right) \right.\\
\left. \phantom{xxxxxx} - \mathbb P_R \left
(  x^\#  + \beta^{(\nell)} K^* ( y - Kx^\#)\right) -  \uell \right
\| =0.
\end{array}
\end{equation}
Denote $h^\#:=x^\#+\beta^{(\infty)} K^* (y-K x^\#)$ and $\hell :=
x^\#  + \beta^{(\nell)} K^* ( y - Kx^\#)$. We have now
\begin{eqnarray*}
\|\P_R(h^\#+\uell)-\P_R(h^\#)-\uell\| \!\!\!\!\!\!\!\!\!\!\!\!\!\!\!\!\!\!\!\!\!\!\!\!\!\!\!\!\!\!\!\!\!\!\!\!\!\!\!\!\!\!\!\!\!\!\!\!\!\!\!\!\!\!&& \\
&&\leq  \|\P_R(\hell+\uell) -\P_R(\hell)-\uell\| \\
&&  + \|\P_R(\hell+\uell)-\P_R(h^\#+\uell)\| + \|\P_R(\hell)-\P_R(h^\#)\|\\
&&\leq \|\P_R(\hell+\uell)-\P_R(\hell)-\uell\|+2\|\hell-h^\#\|.
\end{eqnarray*}
Since both terms on the right hand side converge to zero for $\ell
\to \infty$ (see \eqref{helpeq}), we have
\begin{equation}\label{help2}
\lim_{\ell \to \infty} \|\P_R(h^\#+\uell)-\P_R(h^\#)-\uell\| =0.
\end{equation}
Without loss of generality we can assume $\|h^\#\|_1 >R$. By Lemma
\ref{l1projectionlemma}  there exists  $\mu > 0$ such that
$\P_R(h^\#)= \S_{\mu}(h^\#)$. Because $|h^\#_\lambda| \rightarrow 0$
as $|\lambda| \rightarrow \infty$, this implies that, for some
finite $K_1>0$, $\left(\P_R(h^\#)\right)_\lambda =0 $ for $|\lambda|
> K_1$. Pick now any $\epsilon > 0$ that satisfies $\epsilon <
\mu/5$. There exists a finite $K_2>0$ so that $\sum_{|\lambda| >
K_2} |h^\#_{\lambda}|^2 < \epsilon ^2$. Set $K_0 := \max(K_1, K_2)$,
and define the vector $\tilde h^\# $ by $\tilde h^\#_\lambda
=h^\#_\lambda$ if $|\lambda| \leq  K_0$, $\tilde h^\#_\lambda=0$ if
$|\lambda| > K_0$.

By the weak convergence of the $\uell$, we can, for this same $K_0$,
determine $L_1>0$ such that, for all $\ell \geq L_1$,
$\sum_{|\lambda| \leq K_0} |\uell_\lambda|^2 \leq \epsilon^2$.
Define new vectors $\tuell$ by $\tuell_\lambda=0$ if $|\lambda| \leq
K_0$, $\tuell_\lambda=\uell_\lambda$ if $|\lambda|
> K_0$.

Because of \eqref{help2}, there exists $~L_2~>0$ such that
$\|\P_R(h^\#+\uell)-\P_R(h^\#) -\uell\| \leq \epsilon$ for $\ell
\geq L_2$. Consider now $\ell \geq L:=\max(L_1,L_2)$. We have
\begin{eqnarray*}
&&\|\P_R(\tilde h^\#+\tuell)-\P_R(\tilde h^\#)-\tuell\|  \\
&\leq&\|\P_R(\tilde h^\#+\tuell)-\P_R(h^\#+\tuell)\|
+ \|\P_R(h^\#+\tuell)-\P_R(h^\#+\uell)\|\\
&+& \, \|\P_R(h^\#+\uell)-\P_R(h^\#) -\uell\|
+ \|\P_R(h^\#) -\P_R(\tilde h^\#)\|+\|\uell-\tuell\| \\
 &\leq& ~ 5 \epsilon ~.
\end{eqnarray*}

On the other hand, Lemma \ref{l1projectionlemma} tells us that there
exists $\sigma_\ell >0$ such that $\P_R(\tilde h^\#+\tuell)=
\S_{\sigma_\ell}(\tilde h^\#+\tuell) = \S_{\sigma_\ell}(\tilde
h^\#)+\S_{\sigma_\ell}(\tuell)$, where we used in the last
equality that $\tilde h^\#_\lambda=0$ for $|\lambda| >K_0$ and
$\tuell_\lambda=0$ for $|\lambda| \leq K_0$. From $\|\S_{\mu}(\tilde
h^\#)\|_1 =R =\|\S_{\sigma_\ell}(\tilde
h^\#)\|_1+\|\S_{\sigma_\ell}(\tuell)\|_1$ we conclude that $\sigma_\ell
\geq \mu$ for all $\ell \geq L$. We then deduce
\begin{eqnarray*}
(5 \epsilon)^2 & \geq & \|\P_R(\tilde h^\#+\tuell)-\P_R(\tilde h^\#)-\tuell\|^2\\
& = &\sum_{|\lambda| \leq K_0} |S_{\sigma_\ell}(\tilde h^\#_\lambda)
- S_{\mu}(\tilde h^\#_\lambda)|^2
~+~ \sum_{|\lambda| >K_0}  |S_{\sigma_\ell}(\tuell_\lambda)-\tuell_\lambda |^2 \\
& \geq& \sum_{|\lambda|>K_0} \left[ \max \left(
|\tuell_\lambda|-\sigma_\ell,0\right)
-  |\tuell_\lambda| \right]^2 \\
& = &\sum_{|\lambda|>K_0} \min\left(
|\tuell_\lambda|,\sigma_\ell\right)^2 ~ \geq ~ \sum_{|\lambda|>K_0}
\min\left( |\tuell_\lambda|,\mu\right)^2~.
\end{eqnarray*}

\noindent Because we picked $\epsilon < \mu/5$, this is possible
only if $|\tuell_\lambda| \leq \mu$ for all $|\lambda| > K_0$, $\ell
\geq L$, and if, in addition,
\begin{equation}
\left[
\sum_{|\lambda|>K_0} |\tuell_\lambda|^2\right]^{1/2} ~\leq ~ 5
\epsilon~, ~~~\mbox{i.e.,}~~ \|\tuell\| \leq 5 \epsilon~.
\end{equation}
It then follows that $\|\uell\| \leq \|\tuell\|+
\left[\sum_{|\lambda|\leq K_0} |\uell_\lambda|^2\right]^{1/2} \leq 6
\epsilon$.

\noindent We have thus obtained what we set out to prove: the
subsequence $\left(x^{n_{j_\ell}}\right)_{\ell \in \N}$ of
$(x^{(n)})_{n\in \N}$ satisfies that, given arbitrary $\epsilon >0$,
there exists $L$ so that, for $\ell > L$,
$\,\|x^{n_{j_\ell}}-x^{\#}\| \leq 6 \epsilon$. \hfill $\Box$

\begin{Remark}
In this proof we have implicitly assumed that $\|h^\#+u\j\|_1>R$.
Given that $\|h^\#\|_1>R$, this assumption can be made without loss
of generality, because it is not possible to have $\|h^\#\|_1>R$ and
$\|h^\#+u\j\|_1<R$ infinitely often, as the following argument
shows. Find $K_0,\, L_0$ such that
$\sum_{|\lambda|<K_0}|h^\#_\lambda|\geq (\|h^\#\|_1+R)/2$ and,
 $\forall \ell\geq L_0$ and $\forall |\lambda|<K_0$:
$|\uell_\lambda|<(K_0^{-1}(\|h^\#\|_1-R)/4$. Then
$\sum_{|\lambda|<K_0}|h^\#_\lambda+\uell_\lambda|\geq
\sum_{|\lambda|<K_0}|h^\#_\lambda|-|\uell_\lambda|\geq
(\|h^\#\|_1+R)/2-
(\|h^\#\|_1-R)/4=R+(\|h^\#\|_1-R)/4>R$. Hence, $\forall \ell>L_0$,
$\|h^\#+\uell\|_1\geq R$.
\end{Remark}

\begin{Remark}At the cost of more technicalities it is possible to show that the
whole subsequence $(x^{(n_j)})_{ j \in \mathbb{N}}$ defined in the
proof of Proposition \ref{weakconv} converges in norm to  $x^\#$,
i.e., that $\lim_{j \to \infty} \| x^{(n_j)}- x^\#\| =0$, without
going to a subsequence $\left(x^{n_{j_\ell}}\right)_{\ell \in \N}$.
\end{Remark}

\noindent The following proposition summarizes in one statement all
the findings of the last two subsections.

\begin{Proposition}[Norm  convergence to minimizing accumulation points]
\label{normconv} Every weak accumulation point $x^\#$ of the
sequence $(x\n)_{n \in \mathbb{N}}$ defined by \eqref{PGiteration}
is a minimizer of $\mathcal{D}$ in $B_R$. Moreover, there exists a
subsequence $(x^{(n_\ell)})_{ \ell \in \mathbb{N}}$ of $(x\n)_{n \in
\mathbb{N}}$ that converges to $x^\#$ in norm.
\end{Proposition}

\subsection{Uniqueness of the accumulation point}

In this subsection we prove that the accumulation point $x^\#$ of
$(x\n)_{n \in \mathbb{N}}$ is unique, so that the entire sequence
$(x\n)_{n \in \mathbb{N}}$  converges to $x^\#$ in norm. (Note that
two sequences $(x\n)_{n \in \mathbb{N}}$ and $({x'}\n)_{n \in \mathbb{N}}$,
both defined by the same recursion, but starting from different initial
points $x^{(0)} \neq {x'}^{(0)}$, can still converge to different limits
$x^\#$ and ${x'}^\#$.)


We start again from the inequality
\begin{equation} \label{maininequ}
\langle  x\n + \beta\n K^* ( y - K x\n) - x\np1 , w - x\np1 \rangle
\leq 0,
\end{equation}
for all $w \in B_R$ and for all $n \in \mathbb{N}$,
and its many consequences. Define $M_R$ to be the set of minimizers
of $\mathcal{D}$ on $B_R$. By Lemma \ref{kern}, $M_R = B_R \cap
(\tilde x + \ker K)$, where $\tilde x$ is an arbitrary minimizer of
$\mathcal{D}$ in $B_R$. By the convention adopted in Remark
\ref{conven},
\begin{equation}
M_R \subset B_R^+ :=\left \{x \in \ell_1(\Lambda); x_\lambda \geq 0
\text{ for all } \lambda \in \Lambda, \text{ and } \sum_{\lambda \in
\Lambda} x_\lambda \leq R \right \}.
\end{equation}
Moreover, for each element $z \in M_R$, $z_\lambda=0$ if $\lambda
\notin \Gamma$ (see Lemma \ref{gammaset}). The set $M_R$ is both
closed and convex. We define the corresponding (nonlinear)
projection operator $\P_{M_R}$ as usual,
\begin{equation}
\P_{M_R}(v):=\arg \, \min \{ \|v- z\|^2\,; \, z \in M_R\,\}\,.
\end{equation}
Because $M_R$ is convex, this projection operator has the following
property:
\begin{equation}\label{charproj}
\forall \tilde x \in M_R\,:~\langle z - \P_{M_R}(z), \tilde x -
\P_{M_R}(z) \rangle \leq 0.
\end{equation}
(The proof is standard, and is essentially given in the proof Lemma
\ref{contractionlemma}, where in fact only the convexity of $B_R$ was used.) For each $n \in \mathbb{N}$, we introduce now $a\n$
and $b\n$ defined by
\begin{equation}
a\n := \P_{M_R}(x\n), \quad b\n = x\n -a\n.
\end{equation}
Specializing equation \eqref{charproj} to $x\n$, we obtain, for all
$\tilde x \in M_R$ and for all $n \in \mathbb N$:
\begin{equation}\label{charproj2}
\langle x\n - a\n, \tilde x - a\n \rangle \leq 0.
\end{equation}
or
\begin{equation}\label{charproj3}
\langle b\n, \tilde x - a\n \rangle \leq 0.
\end{equation}
Because $a\n$ is a minimizer, we can also apply Lemma \ref{fixpt} to
$a\n$ and conclude
\begin{equation}\label{charproj4}
\langle K^*(y - K a\n), w - a\n \rangle \leq 0, \quad \text{ for all
} w \in B_R.
\end{equation}

With these inequalities, we can prove the following crucial result.

\begin{Lemma}\label{monolem}
For any $\tilde x \in M_R$, and for any $n \in \mathbb{N}$,
\begin{equation} \label{monot}
\| x\np1 -\tilde x\| \leq \| x\n -\tilde x\|.
\end{equation}
\end{Lemma}
{\em Proof:} We set $w= \tilde x$ in \eqref{maininequ}, leading to
\begin{equation}
\langle x\n - x\np1 , \tilde  x - x\np1 \rangle + \beta\n \langle
K^*(y - K x\n ), - b\np1\rangle \leq 0,
\end{equation}
where we have used that $K \tilde x = K a\np1$. We also have,
setting $w = x\np1$ in the $(n+1)$-version of \eqref{charproj4},
\begin{equation}
\langle K^*(y - K a\np1),x\np1 - a\np1 \rangle \leq 0,
\end{equation}
or
\begin{equation}
\langle K^*(y - K a\n),b\np1 \rangle \leq 0,
\end{equation}
where we have used $K a\n = K a\np1$. It follows that
\begin{equation}
\langle x\n - x\np1 , \tilde  x - x\np1 \rangle + \beta\n \langle
-K^*K b\n, - b\np1\rangle \leq 0,
\end{equation}
or
\begin{equation}
\langle x\n - x\np1 , \tilde  x - x\np1 \rangle + \beta\n \langle K
b\n, K b\np1\rangle \leq 0,
\end{equation}
which is also equivalent to
$$
\langle x\n - \tilde x , \tilde  x - x\np1 \rangle +\|
\tilde x - x\np1\|^2 + \frac{1}{2}\beta\n \left [ \| K
b\n\|^2 + \| K b\np1\|^2 \right
]$$
\begin{equation}
\label{addto}
-\frac{1}{2} \beta\n \| K b\n - K b\np1\|^2 \leq 0.
\end{equation}
Adding $\frac{1}{2} \beta\n \| K(b\n - b\np1)\|^2 \leq
\frac{r}{2} \| x\n - x\np1\|^2$ to \eqref{addto}, we have
\begin{eqnarray*}
&& \langle x\n - \tilde x , \tilde  x - x\np1 \rangle +\| \tilde x -
x\np1\|^2 + \frac{1}{2}\beta\n \left [ \| K b\n\|^2 + \| K b\np1\|^2\right ]\\
 &\leq& \frac{r}{2} \| x\n - x\np1\|^2\\
&=& \frac{r}{2} \left [ \| x\n- \tilde x \|^2 + \| x\np1 - \tilde x
\|^2 - 2 \langle  x\n -\tilde x, x\np1 - \tilde x \rangle \right ].
\end{eqnarray*}
It follows that
$$
\left ( 1 - \frac{r}{2} \right ) \| x\np1- \tilde x \|^2 + (1 -r)
\langle \tilde x - x\n, x\np1 - \tilde x \rangle -\frac{r}{2} \|
\tilde x - x\n \|^2$$
\begin{equation}
\leq - \frac{1}{2}  \beta\n \left [ \| K
b\n\|^2 + \|K b\np1\|^2 \right ] \leq 0,
\end{equation}
which, in turn, implies that
\begin{equation}
\left ( 1 - \frac{r}{2} \right )  \| x\np1- \tilde x \|^2 - (1 -r)\|
\tilde x - x\n\| \|x\np1 - \tilde x\| - \frac{r}{2} \| \tilde x -
x\n \|^2 \leq 0.
\end{equation}
This can be rewritten as
\begin{equation}
\left [ \| \tilde x - x\np1\| - \|\tilde x - x\n\| \right ]  \left[
 \left ( 1 - \frac{r}{2} \right ) \| x\np1- \tilde x \|
+ \frac{r}{2}  \|\tilde x - x\n\|  \right ] \leq 0,
\end{equation}
which implies $\| x\np1 -\tilde x\| \leq \| x\n -\tilde x\|$.
 \hfill $\Box$

\noindent We are now ready to state the main result of our work.
\begin{Theorem}
The sequence $\left(x\n\right)_{n \in \mathbb{N}}$ as defined in
\eqref{PGiteration}, where the {\em step-length sequence}
$\left(\beta\n\right)_{n \in \mathbb{N}}$ satisfies Condition (B)
with respect to the $x\n$, converges in norm to a minimizer of
$\mathcal{D}$ on $B_R$. \label{finaltheorem}
\end{Theorem}
{\em Proof:} The sequence $(x\n)_{n \in \mathbb{N}}$ has a least one
accumulation point $x^\#$. By Proposition \ref{weakconv}  $x^\#$
minimizes $\mathcal{D}$ in $B_R$. By Proposition \ref{normconv}
$(x\n)_{n \in \mathbb{N}}$ has a subsequence
$\left(x^{(n_\ell)}\right)_{ \ell \in \mathbb{N}}$ that converges to
$x^\#$. By Lemma \ref{monolem} $\|x\n - x^\#\|$ decreases
monotonically, hence it has a limit for $n \to \infty$, and
\begin{equation}
\lim_{n \to \infty } \|x\n - x^\#\| = \lim_{\ell \to \infty}
\|x^{(n_\ell)} - x^\#\| = 0.
\end{equation}
\hfill $\Box$

\section{Numerical Experiments and Additional Algorithms}

\label{numericsection}

\subsection{Numerical examples}

We conduct a number of numerical experiments to gauge the
effectiveness of the different algorithms we discussed.
All computations were done in Mathematica 5.2
\cite{mathematica5} on a 2Ghz workstation with 2Gb memory.

We are
primarily interested in the behavior, as a function of time (not
number of iterations), of the relative error $\|x^{(n)}-\bar
x\|/\|\bar x\|$. To this end, and for a given operator $K$ and data
$y$, we need to know in advance the actual minimizer $\bar
x(\tau)$ of the functional (\ref{functional}).


One can calculate the minimizer exactly (in practice up to computer
round-off) with a finite number of steps using the LARS algorithm
described in \cite{efhajoti04} (the variant called `Lasso',
implemented independently by us). This algorithm scales badly, and
is useful in practice only when the number of non-zero entries in
the minimizer $\bar x(\tau)$ is sufficiently small. We made our own
implementation of this algorithm to make it more directly applicable
to our problem (i.e., we do not renormalize the columns of the
matrix to have zero mean and unit variance, as it is done in the
statistics context \cite{efhajoti04}). We also double-check the
minimizer obtained in this manner by verifying that it is indeed a
fixed point of the iterative thresholding algorithm
(\ref{DDD04iteration}) (up to machine epsilon). We then have an
`exact' minimizer $\bar x$ together with its radius $R=\|\bar x\|_1$
(used in the projected algorithms) and, according to Lemma
\ref{thrval}, the corresponding threshold $\tau=\max_i|\bar r_i|$
with $\bar r=K^\ast(y-K\bar x)$ (used in the
iterative thresholding algorithm).\\

The numerical examples below are listed in order of
increasing complexity; they illustrate that the algorithms can behave
differently for different examples. In these experiments we choose
$\beta^{(n)}= \beta\n_{\mbox{\tiny{st.}}} := \|r^{(n)}\|^2/\|Kr^{(n)}\|^2$,
(where, as before,  $r^{(n)}=K^\ast
(y-K x\n)$); $\beta\n_{\mbox{\tiny{st.}}}$ is the
standard descent parameter from the classical linear steepest
descent algorithm.
\begin{enumerate}
\item When $K$ is a partial Fourier matrix (i.e., a Fourier
matrix with a prescribed number of deleted rows), there is no
advantage in using a dynamical step size
$\beta_{\mbox{\tiny{st.}}}^{(n)}=\|r^{(n)}\|^2/\|Kr^{(n)}\|^2$ as this ratio is always
equal to 1. This trivially fulfills Condition (B) in Section
\ref{genpropsec}.  The performance of the projected steepest descent
iteration simply equals that of the projected Landweber
iterations.
\item By combining a scaled partial Fourier transform with a rank 1
projection operator, we constructed our second example, in which
$K$ is a $1536 \times 2049$ matrix, of rank $1536$, with largest
singular value
equal to 0.99 and all the other singular values between 0.01 and
0.11. Because of the construction of the matrix,
the FFT algorithm provides a fast way of computing
the action of this matrix on a vector.
For the $y$ and $\tau$ that were chosen,
the limit vector $\bar{x}_{\tau}$ has $429$ nonzero entries.
For this example, the LARS
procedure is slower than thresholded Landweber, which in turn is
significantly slower than projected steepest descent.
To get within a distance of the true minimizer corresponding to a $5\%$ relative error,
the projected steepest descent algorithm
takes $2\sec$, the thresholded Landweber algorithm $39\sec$, and LARS
$151\sec$. (The relatively poor performance of LARS in this case is
due to the large number of nonzero entries in the limit vector
$\bar{x}_{\tau}$; the complexity of LARS is cubic in this number
of nonzero entries.) In this case, the
$\beta^{(n)}_{\mbox{\tiny{st.}}}=\|r^{(n)}\|^2/\|Kr^{(n)}\|^2$ are much larger
than 1; moreover, they satisfy Condition (B)
of Section \ref{genpropsec} at every step.  We illustrate the
results in Figure \ref{PSDfig}.
\begin{figure}[ht]
\resizebox{\textwidth}{!}{\includegraphics{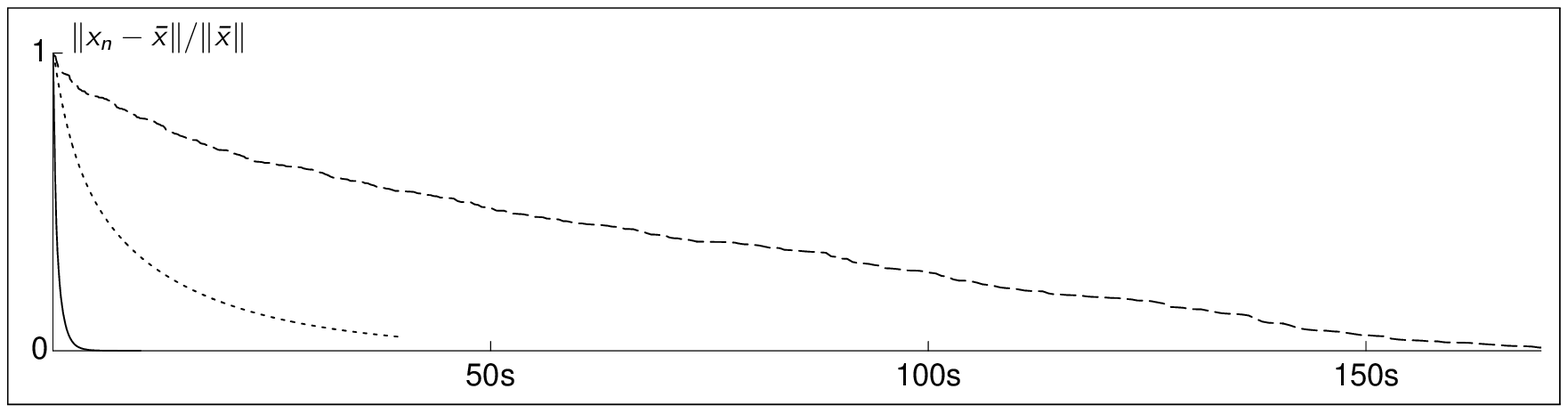}}
\caption{The different convergence rates of the thresholded
Landweber algorithm (dotted line), the projected steepest descent
algorithm (solid line, near vertical axis) and the LARS algorithm
(dashed line), for the second example.
The projected steepest descent algorithm converges
much faster than the thresholded Landweber iteration. They both do
better than the LARS method.}\label{PSDfig}
\end{figure}
\item The last example is inspired by a real-life application in
geoscience \cite{LoNoDaDa}, in particular an application in seismic
tomography based on earthquake data. The object space consists of the
wavelet coefficients of a 2D seismic velocity perturbation. There
are $8192$ degrees of freedom. In this particular case the number of
data is $1848$. Hence the matrix $K$ has $1848$ rows and $8192$
columns. We apply the different methods to the same noisy data that
are used in \cite{LoNoDaDa} and measure the time 
to convergence up to a specified relative error (see Table
\ref{geodistancetable} and Figure \ref{geodistancefig}). This example
illustrates the slow convergence of the thresholded Landweber
algorithm (\ref{DDD04iteration}), and the improvements made by a
projected steepest descent iteration (\ref{PGiteration}) with
the  special choice $\beta^{(n)}=\beta\n_{\mbox{\tiny{st.}}}$ above.
In this case, this choice turns out {\em not} to satisfy Condition (B) in
general. One could conceivably use successive corrections, e.g. by a line-search,
to determine,
starting from $\beta\n_{\mbox{\tiny{st.}}}$, values of $\beta\n $ that would
satisfy condition (B), and thus
guarantee convergence as established by Theorem 5.18.
This would slow down the method considerably. The $\beta\n_{\mbox{\tiny{st.}}}$
seem to be in the right ballpark, and provide us with a numerically converging sequence.
We also implemented the projected Landweber algorithm
(\ref{PLiteration}); it is listed in Table
\ref{geodistancetable} and illustrated in Figure \ref{geodistancefig}.

The matrix $K$ in this example is extremely
ill-conditioned: its largest singular value was normalized to 1,
but the remaining singular values quickly tend to zero. The
threshold was chosen, according to the (known or estimated)
noise level in the data,
so that $\mathcal{D}(\bar
x)/\sigma^2=1848$ ( = the number of data points),
where $\sigma$ is the data noise level; this is a standard
choice that avoids overfitting.

\begin{figure}
\resizebox{\textwidth}{!}{\includegraphics{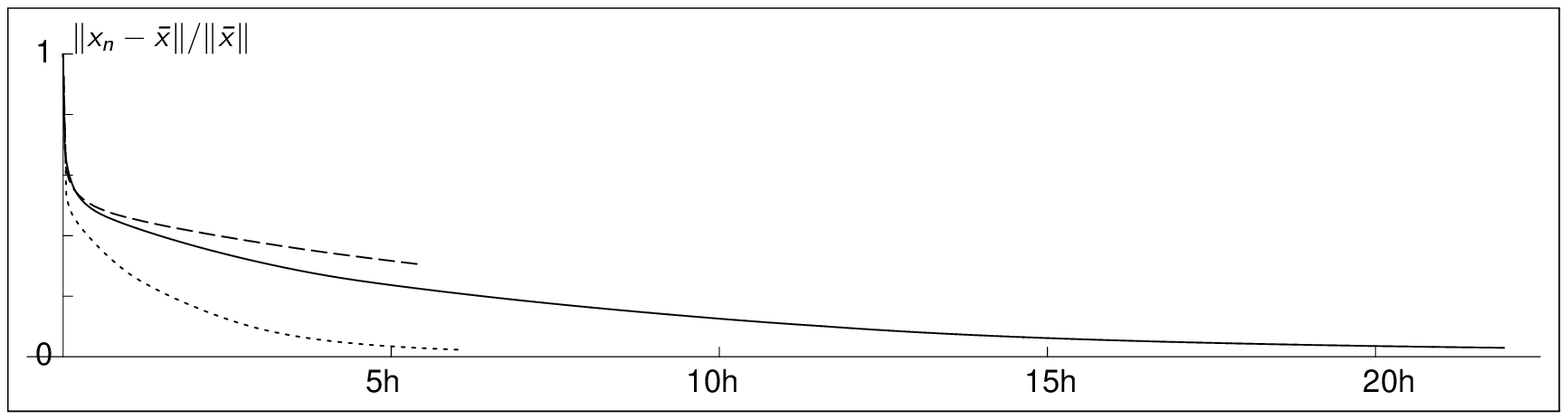}}
\caption{The different convergence rates of the thresholded
Landweber algorithm (solid line), the projected Landweber algorithm
(dashed line) and the projected steepest descent algorithm (dotted
line), for the third example.
The projected steepest descent algorithm converges about four
times faster than the thresholded Landweber iteration. The projected
Landweber iteration does better at first (not visible in this plot),
but looses with respect to iterative thresholding afterwards. The
horizontal axis has time (in hours), the vertical axis displays the
relative error.}\label{geodistancefig}
\end{figure}

In Figure \ref{geodistancefig}, we see that the thresholded
Landweber algorithm takes more than $21$ hours (corresponding to
$200,000$ iterations) to converge to the true minimizer within a $3\%$ relative error,
as measured by the usual $\ell_2$ distance.
The projected steepest descent algorithm is about four times faster
and reaches the same reconstruction error in about $5.5$ hours
($25,000$ iterations). Due to one additional matrix-vector
multiplication and, to a minor extent, the computation of the
projection onto an $\ell_1$-ball, one step in the projected steepest
descent algorithm takes approximately twice as long as one step in
the thresholded Landweber algorithm. For the projected Landweber
algorithm there is an advantage in the first few iterations,
but after a short while,
the additional time needed to compute the projection $\P_R$ (i.e.,
to compute the corresponding variable thresholds) makes this
algorithm slower than the iterative soft-thresholding. We illustrate
the corresponding CPU time in Table \ref{geodistancetable}.

It is worthwhile noticing that for the three algorithms the value of the
functional (\ref{functional}) converges much faster to its limit
value than the minimizer itself: When the reconstruction error
is 10\%, the corresponding value of the functional is already
accurate up to three digits with respect to the value of the
functional at $\bar x$. We can imagine that in this case the
functional has a long narrow ``valley'' with a very gentle slope in
the direction of the eigenvectors with small (or zero) singular
values.
\begin{table}
\centering\begin{tabular}{r|rr|rr|rr}
Relative & \multicolumn{2}{c|}{thresholded Landweber} &
\multicolumn{2}{c|}{projected st. descent} &
\multicolumn{2}{c}{projected Landweber}\\
error & $n$ & time & $n$ & time & $n$ & time\\
\hline 0.90 & 3 & 1s & 2 & 1s & 3 & 2s \\
0.80 & 20 & 8s & 8 & 7s &  15 & 11s \\
0.70 & 163 & 1m8s & 20 & 17s & 59 & 44s \\
0.50 & 3216  & 22m9s & 340 & 4m56s & 2124 & 27m17s \\
0.20 & 55473 & 6h23m & 6738 & 1h37m &  &  \\
0.10 & 100620 & 11h38m & 11830 & 2h51m &  &\\
0.03 & 198357 & 21h47m & 22037 & 5h20m &  &
\end{tabular}
\caption{Table illustrating the relative performance of three
algorithms: thresholded Landweber, projected Landweber and projected
steepest descent, for the third example.}\label{geodistancetable}
\end{table}
%
The path in the  $\|x\|_1$ vs. $\|Kx-y\|^2$ plane followed by the
iterates is shown in Figure \ref{pathsfig1}. The projected steepest
descent algorithm, by construction, stays within a fixed
$\ell_1$-ball, and, as already mentioned, converges faster than the
thresholded Landweber algorithm. The path followed by the LARS
algorithm is also pictured. It corresponds with the so-called {\it
trade-off} curve which can be interpreted as the border of the area
that is reachable by any element of the model space, i.e., it is
generated by $\bar x(\tau)$ for decreasing values of $\tau>0$.

In this particular example, the number of nonzero components of
$\bar x$ equals $128$. The LARS (exact) algorithm only takes $55$
seconds, which is \emph{much} faster than any of the iterative
methods demonstrated here. However, as illustrated above, by the second example,
LARS looses its advantage
when dealing with larger problems where the minimizer is not sparse in absolute number of entries, as is the case
in, e.g., realistic problems of global seismic tomography. Indeed, the
example presented here is a ``toy model'' for proof-of-concept for
geoscience applications. The 3D model will involve millions of
unknowns and solutions that may be sparse compared with the total
number of unknowns, but not sparse in absolute numbers. Because the complexity
of LARS is cubic
in the number of nonzero components of the
solution, such 3D problems are expected to lie beyond its useful range.
\end{enumerate}


\subsection{Relationship to other methods}

The projected iterations (\ref{projLWiteration}) and
(\ref{projSDiteration}) are 
related to the POCS  (Projection
on Convex Sets) technique \cite{Br65}. The projection of a vector
$a$ on the solution space $\{x: Kx=y\}$ (a convex set, assumed here to
be non-empty; no such assumption was made before because the
functional (\ref{functional}) always has a minimum) is given by:
\begin{equation}
x=a-K^* (KK^*)^{-1}(y-Ka) \label{projection}
\end{equation}
Hence, alternating projections on the convex sets  $\{x: Kx=y\}$ and
$B_R$ give rise to the algorithm \cite{CaRo04}: : Pick an arbitrary $x^{(0)} \in \ell_2(\Lambda)$, for example $x^{(0)}=0$, and iterate
\begin{equation}
x\np1=\P_R(x\n-K^* (KK^*)^{-1}(y-Kx\n))
\end{equation}
This may be practical in case of a small number of data or when
there is structure in $K$, i.e., when $KK^*$ is efficiently inverted.
Approximating $KK^\ast$ by the unit matrix, yields the projected
Landweber algorithm (\ref{projLWiteration}); approximating
$(KK^*)^{-1}$ by a constant multiple of the unit matrix yields the
projected gradient iteration (\ref{projSDiteration}) if one chooses
the constant equal to $\beta\n$.
\begin{figure}[ht]
\begin{center}
\resizebox{10cm}{!}{\includegraphics{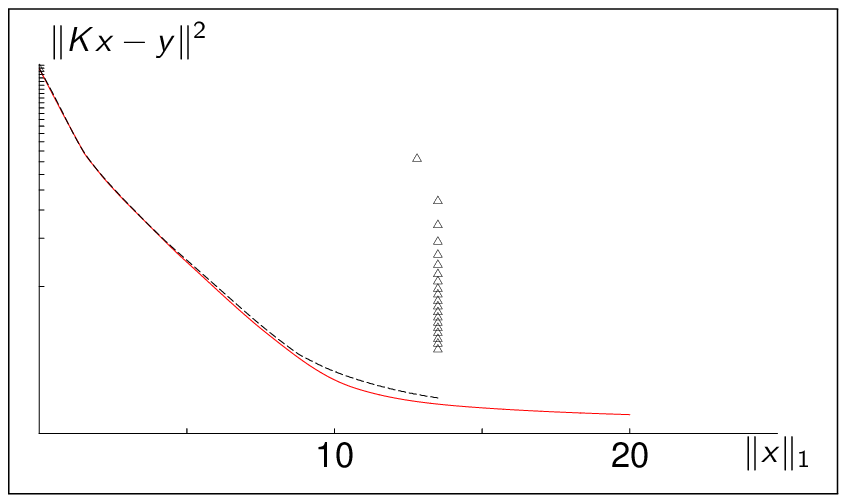}}
\caption{Trade-off curve (solid line) and its approximation with algorithm
\eqref{adaptiveSDiteration} in 200 steps (dashed line). For comparison, the iterates
of projected steepest descent are also indicated (triangles).}
\end{center}
\label{adaptivefigure}
\end{figure}
The 
projected methods discussed in this paper produce iterates
that (except for the first few) live on the `skin' of the
$\ell_1$-ball of radius $R$, as shown in Fig. \ref{pathsfig1}.
We have found even more promising results for an  `interior' algorithm
in which we still project on $\ell_1$-balls, but now with a slowly
increasing radius, i.e.,
\begin{equation}
x\np1=\P_{R\n}\left(x\n+\beta^{(n)}r\n\right), \quad R\n=(n+1)R/N, \mbox{ and } n=0,\dots, N,
\label{adaptiveSDiteration}
\end{equation}
where $N$ is the prescribed maximum number of iterations (the
origin is chosen as the starting point of this iteration). We do not
have a proof of convergence of this `interior point type'
algorithm. 
We observed (also without proof) that the
path traced by the iterates $x\n$ (in the
 $\|x\|_1$ vs. $\|Kx-y\|^2$ plane) 
is very close to the trade-off curve (see Fig. 5);
this is a useful property in practice since at least part of the trade-off
curve should be constructed anyway.\\
Note that the strategy followed by these algorithms is similar to that of LARS
\cite{efhajoti04}, in that they both start with $x^{(0)}=0$ and slowly
increase the $\ell_1$ norm of the successive approximations. It is also related to \cite{Hale.Yin.ea2007}.\\
While we were finishing this paper, Michael Friedlander informed us
of their numerical results in \cite{VF} which are closely related to
our approach,
although their analysis is limited to finite dimensions.\\
Different, but closely related is also the recent approach by
Figueiredo, Nowak, and Wright \cite{finowr07}. The authors first
reformulate the minimization of \eqref{functional} as a
bound-constrained quadratic program in standard form, and then they
apply iterative projected gradient iterations, where the projection
act componentwise by clipping to zero negative components.

\section{Conclusions}

We have presented  convergence results for accelerated projected
gradient methods to find a minimizer of an $\ell_1$ penalized
functional. The innovation due to the introduction of `Condition
(B)' is to guarantee strong convergence for the full sequence.
Numerical examples confirm that this algorithm can
outperform (in terms of CPU time) existing methods such as the
thresholded Landweber iteration or even LARS.

It is important to remark that the speed of convergence may
depend strongly on how the operator is available.
Because most of the time in the
iterations is consumed by matrix-vector multiplications
(as is often the case for iterative algorithms), it makes a
big difference whether
$K$ is given by a
full matrix or a sparse matrix (perhaps sparse in the sense that
its action on a vector can be computed via a fast algorithm, such
as the FFT or a wavelet transform).
The applicability of the projected algorithms hinges on
the observation that the $\ell_2$ projection on an $\ell_1$ ball can
be computed with a $\mathcal{O}(m \log m)$-algorithm, where $m$ is the dimension
of the underlying space.

There is no universal method that
performs best for any choice of the operator, data, and penalization
parameter.
As a general rule of thumb we expect that,
among the algorithms discussed in this paper
for which we have convergence proofs,
\begin{itemize}
\item the thresholded Landweber algorithm
(\ref{DDD04iteration}) works best for an operator $K$ close to the
identity (independently of the sparsity of the limit),
\item the projected steepest descent algorithm \eqref{PGiteration} works
best for an operator with a relatively nice spectrum, i.e., with not
too many zeroes (also independently of the sparsity of the
minimizer), and
\item the exact (LARS) method works best when the minimizer is
sparse in absolute terms.
\end{itemize}
Obviously, the three cases overlap partially,
and they do not cover the whole range of possible operators and
data. In future work we intend to investigate algorithms that
would further improve the performance for the case of a large ill-conditioned
matrix and a minimizer that is relatively sparse with respect to the
dimension of the underlying space. We intend, in particular, to focus
on proving convergence and other mathematical properties of
\eqref{adaptiveSDiteration}.

\section{Acknowledgments}

M.~F. acknowledges the financial support provided by the European
Union's Human Potential Programme under the contract
MOIF-CT-2006-039438. I.~L. is a post-doctoral fellow with the
F.W.O.-Vlaanderen (Belgium). M.F. and I.L. thank the Program in
Applied and Computational Mathematics, Princeton University, for the
hospitality during the preparation of this work. I.~D. gratefully
acknowledges partial support from NSF grants DMS-0245566 and
0530865.


\providecommand{\bysame}{\leavevmode\hbox to3em{\hrulefill}\thinspace}
\providecommand{\MR}{\relax\ifhmode\unskip\space\fi MR }
\providecommand{\MRhref}[2]{%
  \href{http://www.ams.org/mathscinet-getitem?mr=#1}{#2}
}
\providecommand{\href}[2]{#2}

\end{document}